\IfFileExists{elsarticle.cls}{%
  \documentclass[review]{elsarticle}%
}{%
  \documentclass[11pt]{article}
  \input{elsarticle-fallback.tex}
}
\usepackage[dvipsnames]{xcolor}
\usepackage{amssymb,amsthm,amsmath}
\usepackage[margin=2.5cm, bottom=2cm]{geometry}
\usepackage{multirow}
\usepackage{subcaption}
\usepackage{booktabs}
\usepackage[colorlinks=true, linkcolor=blue, citecolor=red, urlcolor=blue]{hyperref}
\usepackage{lineno}
\journal{}

\newcommand{\Eqref}[1]{Equation~\eqref{#1}}
\newcommand{\Figref}[1]{Figure~\ref{#1}}

\newcommand{\Sref}[1]{Section~\ref{#1}}
\newcommand{\Aref}[1]{\ref{#1}}

\begin{document}
\begin{frontmatter}

\title{Tensor-Train Isogeometric Analysis for Large-Scale Three-Dimensional Linear Elasticity}

\author[T]{Quoc Thai Tran\texorpdfstring{\corref{cor1}}{}}
\author[T]{Duc P. Truong}
\author[XCP]{William Dai}
\author[T]{Kim \O. Rasmussen}
\address[T]{Theoretical Division, Los Alamos National Laboratory, Los Alamos, NM, USA}
\address[XCP]{X Computational Physics Division, Los Alamos National Laboratory, Los Alamos, NM, USA}
\cortext[cor1]{Corresponding author\\
Email addresses: thai.tran@lanl.gov\\LA-UR-26-25203}
\begin{abstract}
This work develops a tensor-train (TT) isogeometric analysis (IGA) framework for large-scale three-dimensional linear elasticity. The tensor-product structure of IGA basis functions is exploited to construct the geometry-dependent elasticity operator, force vector, and displacement solution directly in compressed TT format. In contrast to conventional full-grid assembly, the proposed formulation avoids storing the complete multidimensional discrete operators and instead represents the spatial directions and displacement-component coupling through low-rank TT cores. Geometry-dependent coefficient fields are approximated by TT-cross interpolation, and the resulting coupled elasticity system is solved directly in TT format using the alternating minimal energy (AMEn) method. A manufactured three-dimensional elasticity problem verifies that the TT formulation preserves the accuracy and convergence behavior of the conventional full-grid IGA discretization. Further numerical examples demonstrate substantial compression of the elasticity operators and displacement fields and enable simulations beyond the problem sizes accessible to the assembled full-grid formulation. The results demonstrate the potential of TT-IGA as a memory-efficient approach for large-scale three-dimensional elasticity problems whose geometry and solution admit low-rank tensor representations.
\end{abstract}

\end{frontmatter}

\section*{Highlights}
\begin{itemize}
    \item A tensor-train formulation is developed for large-scale three-dimensional elasticity.
    \item Elasticity operators are built directly from the IGA-based geometry in compressed tensor form.
    \item The method preserves full-grid accuracy and enables simulations with up to 1.6 billion DOFs.
\end{itemize}


\section{Introduction}

Three-dimensional linear elasticity is a fundamental model in computational solid mechanics for predicting deformation, stress, and strain in engineering structures and materials \cite{belytschko2014nonlinear,sadd2009elasticity}. At fine spatial resolutions, however, conventional discretizations lead to large algebraic systems whose assembly, storage, and solution become increasingly expensive. This challenge is particularly pronounced in three dimensions because the number of degrees of freedom grows rapidly with refinement.

Isogeometric analysis (IGA) provides an attractive framework for high-fidelity structural analysis by using spline basis functions originating from computer-aided design for both geometry representation and numerical approximation \cite{hughes2005isogeometric,cottrell2009isogeometric}. Its use of B-splines and non-uniform rational B-splines (NURBS) enables accurate geometric descriptions together with smooth high-order approximation spaces. IGA has been successfully extended to complex engineering applications, including patient-specific NURBS modeling for vascular flow \cite{zhang2007patient}, coupled fluid--structure interaction \cite{bazilevs2008isogeometric, bazilevs2006isogeometric}, immersed analysis of complex three-dimensional CAD geometries \cite{antolin2022quadrature}, and CAD-to-analysis pipelines for generating spline representations suitable for numerical simulation \cite{marsala2024cad}. More recently, M\"oller et al.\ introduced IGANets, combining spline-based geometric representations with physics-informed learning for structural analysis \cite{moller2026iganets}. These developments illustrate the continuing effort to preserve accurate geometric representations while improving the flexibility and efficiency of computational analysis.

Despite these advantages, the tensor-product structure underlying standard spline discretizations can lead to rapidly increasing storage and computational requirements in three-dimensional simulations. Low-rank tensor representations provide one possible route for alleviating this difficulty by exploiting separability in the geometry, operators, and solution fields. Among these representations, the tensor-train (TT) format \cite{oseledets2011tensor} expresses a multidimensional tensor as a sequence of low-dimensional cores, with storage governed primarily by the TT ranks rather than by the total number of full-grid entries. When the underlying quantities admit sufficiently low TT ranks, substantial reductions in memory and computational cost can be achieved. TT methods have been applied to a broad range of scientific and engineering problems, including Poisson equations \cite{tran2026tensor}, neutron transport \cite{truong2024tensor,ortega2026Applied}, Maxwell and Vlasov--Maxwell systems \cite{manzini2023tensor,chinomona2026VM}, the Fokker--Planck equation \cite{dolgov2012fast}, elliptic PDEs \cite{kazeev2013low}, stochastic finite-volume methods \cite{walton2025tensor}, chemical master equations \cite{vo2017adaptive}, compressible flows \cite{danis2025tensor}, quantum dynamics \cite{gelss2025quantum}, configurational integrals \cite{truong2025breaking}, and large-scale data analysis \cite{bhattarai2020distributed}. A broad review of low-rank tensor methods for PDEs can be found in \cite{bachmayr2023low}.

The tensor-product structure of IGA is particularly compatible with such low-rank approximations. Ion et al.\ developed a TT-based IGA formulation for PDEs on parameter-dependent geometries \cite{ion2022tensor}, Markeeva et al.\ proposed a quantized TT isogeometric solver \cite{markeeva2021qtt}, and Montardini et al.\ developed a Tucker-based low-rank IGA solver for Poisson and elasticity problems \cite{montardini2023low}. More recently, Riemer and Stoll introduced a TT framework for hierarchical IGA \cite{riemer2026low}. Tran et al.\ \cite{tran2026tensor} recently developed a fully TT-assembled three-dimensional IGA formulation for the Poisson equation, in which geometry-dependent coefficients and discrete operators are constructed directly in TT format while retaining the original NURBS geometry representation.

Low-rank tensor formulations for elasticity remain comparatively less developed. Benvenuti et al.\ recently proposed a quantized tensor-train finite element method for two-dimensional linear elasticity \cite{benvenuti2025low}. Their approach reorganizes the finite element discretization and degree-of-freedom ordering to preserve the QTT structure. Extending a direct TT-IGA formulation from a scalar PDE to three-dimensional elasticity presents a different challenge: the displacement field contains three coupled components, the constitutive relation generates multiple geometry-dependent coefficient fields, and the resulting operator couples derivatives in different parametric directions. These features require the scalar TT operator construction to be generalized to a fully coupled vector-valued system while preserving a compact tensor representation.

To the best of our knowledge, a direct full-order TT-IGA formulation for coupled three-dimensional linear elasticity, in which the elasticity operator is constructed directly from the spline geometry and solved as a coupled TT system, has not previously been reported. The novelty of the present work lies in extending direct TT-IGA assembly from scalar operators to the fully coupled three-dimensional elasticity operator while retaining the original tensor-product IGA indexing and NURBS geometry representation. Geometry-dependent coefficient fields are approximated through TT-cross interpolation and used to construct the one-dimensional operator cores, while the displacement components are incorporated as an additional TT mode. The resulting low-rank system is solved using the alternating minimal energy (AMEn) method. The formulation is assessed against conventional full-grid IGA in terms of numerical accuracy, storage, and computational cost. A manufactured problem verifies the convergence and accuracy of the formulation, while open-hemisphere and perforated-plate examples demonstrate its compression and scalability. The method targets tensor-product parameterizations for which the geometry-dependent fields and displacement solution admit sufficiently low TT ranks; within this class of problems, it substantially extends the problem sizes accessible to IGA without sacrificing the accuracy of the underlying discretization.

\section{Notation and tensor-train preliminaries}
\label{sec.NotationTT}

The notation used throughout this paper is as follows. Vectors are denoted by bold lowercase letters, such as \(\boldsymbol{u}\), and matrices are denoted by bold uppercase letters, such as \(\boldsymbol{A}\). Their corresponding tensor-train representations are denoted by \(\boldsymbol{u}^{\mathtt{TT}}\) and \(\boldsymbol{A}^{\mathtt{TT}}\), respectively. When referring to individual scalar entries, the boldface is omitted. For example, \(\boldsymbol{A}\) denotes a matrix, \(A_{ij}\) denotes its \((i,j)\)-th entry, and
\[
\left(A^{\mathtt{TT}}\right)_{i_1j_1\cdots i_dj_d}
\]
denotes an entry of its tensorized TT representation. Generic higher-order tensors may be denoted by calligraphic uppercase letters, such as \(\mathcal{A}\). TT cores are denoted by symbols such as \(G^{(k)}\), \(Q^{(k)}\), or \(T^{(k)}\), where the superscript \(k\) identifies the corresponding tensor mode. For completeness, consider a \(d\)-dimensional vector tensorized into mode sizes \(n_1,\ldots,n_d\). Its TT representation is written as
\begin{equation}
\left(u^{\mathtt{TT}}\right)_{i_1\cdots i_d}
=
G^{(1)}_{1i_1\alpha_1}
G^{(2)}_{\alpha_1 i_2\alpha_2}
\cdots
G^{(d)}_{\alpha_{d-1}i_d1},
\end{equation}
where
\[
G^{(k)}
\in
\mathbb{R}^{r_{k-1}\times n_k\times r_k},
\qquad
r_0=r_d=1.
\]
The integers \(r_1,\ldots,r_{d-1}\) are referred to as the TT ranks. They determine the storage and computational complexity of the TT representation. The number of scalar entries stored in a TT vector is
\begin{equation}
s_{\boldsymbol{u}^{\mathtt{TT}}}
=
\sum_{k=1}^{d}
r_{k-1}n_kr_k.
\end{equation}
Similarly, for a matrix
\[
\boldsymbol{A}
\in
\mathbb{R}^{(n_1\cdots n_d)\times(m_1\cdots m_d)},
\]
the corresponding TT-matrix representation is written as
\begin{equation}
\left(A^{\mathtt{TT}}\right)_{i_1j_1\cdots i_dj_d}
=
G^{(1)}_{1i_1j_1\alpha_1}
G^{(2)}_{\alpha_1i_2j_2\alpha_2}
\cdots
G^{(d)}_{\alpha_{d-1}i_dj_d1},
\end{equation}
where
\[
G^{(k)}
\in
\mathbb{R}^{r_{k-1}\times n_k\times m_k\times r_k}.
\]
The storage required by the TT-matrix is therefore
\begin{equation}
s_{\boldsymbol{A}^{\mathtt{TT}}}
=
\sum_{k=1}^{d}
r_{k-1}n_km_kr_k.
\end{equation}
Thus, when the TT ranks remain small relative to the full tensor dimensions, TT representations can require substantially fewer stored entries than their corresponding full-grid representations. The storage measures and compression ratios used in the numerical examples are defined in detail in \Aref{sec.Appen-Storage}. The Einstein summation convention is adopted throughout the formulation: repeated indices appearing twice in a tensor expression imply summation over their range, unless otherwise stated. The symbol \(\otimes\) denotes the Kronecker product. For
\[
\boldsymbol{A}\in\mathbb{R}^{n_1\times m_1},
\qquad
\boldsymbol{B}\in\mathbb{R}^{n_2\times m_2},
\]
their Kronecker product satisfies
\[
\boldsymbol{A}\otimes\boldsymbol{B}
\in
\mathbb{R}^{(n_1n_2)\times(m_1m_2)}.
\]
Similarly, if
\[
\boldsymbol{C}\in\mathbb{R}^{n_3\times m_3},
\]
then
\[
\boldsymbol{A}\otimes\boldsymbol{B}\otimes\boldsymbol{C}
\in
\mathbb{R}^{(n_1n_2n_3)\times(m_1m_2m_3)}.
\]
The same notation is used when spatial TT operators are extended by additional displacement-component modes, with the ordering of the resulting tensor specified explicitly in each case.

\section{Theory of linear elasticity}
Let \(\Omega \subset \mathbb{R}^{3}\) denote the physical domain occupied by the elastic body. The displacement field is $\boldsymbol{u}(\boldsymbol{x})$, where $\boldsymbol{x} = \{x_1, x_2, x_3\}$. Under the small-strain assumption, the infinitesimal strain tensor is defined as the symmetric part of the displacement gradient,
\begin{equation}
\boldsymbol{\varepsilon}(\boldsymbol{u}) = \frac{1}{2}\left(\nabla \boldsymbol{u} + \nabla \boldsymbol{u}^{T}\right).
\label{eq.Kinematics}
\end{equation}
The stress tensor is related to the strain tensor through the fourth-order elasticity tensor \(\mathbb{C}\),
\begin{equation}
\boldsymbol{\sigma}(\boldsymbol{u}) = \mathbb{C}:\boldsymbol{\varepsilon}(\boldsymbol{u}).
\end{equation}
For an isotropic homogeneous material, Hooke's law can be written in terms of the Lamé parameters \(\lambda\) and \(\mu\) as
\begin{equation}
\boldsymbol{\sigma}(\boldsymbol{u}) =
\lambda \operatorname{tr}\!\left(\boldsymbol{\varepsilon}(\boldsymbol{u})\right)\boldsymbol{I}
+ 2\mu \boldsymbol{\varepsilon}(\boldsymbol{u}),
\label{eq.Constitutive}
\end{equation}
where \(\boldsymbol{I}\) is the identity tensor. If the material is parameterized by Young's modulus \(E\) and Poisson's ratio \(\nu\), then
\begin{equation}
\mu=\frac{E}{2(1+\nu)},\qquad
\lambda=\frac{E\nu}{(1+\nu)(1-2\nu)}.
\end{equation}
The strong form of the static linear elasticity problem is obtained from the balance of linear momentum. Given a body force \(\boldsymbol{f}\), a prescribed displacement \(\bar{\boldsymbol{u}}\) on the Dirichlet boundary \(\Gamma_D\), and a prescribed traction \(\bar{\boldsymbol{t}}\) on the Neumann boundary \(\Gamma_N\), the governing equations are
\begin{equation}
\begin{aligned}
-\nabla \cdot \boldsymbol{\sigma}(\boldsymbol{u}) &= \boldsymbol{f}
&&\text{in } \Omega,\\
\boldsymbol{u} &= \bar{\boldsymbol{u}}
&&\text{on } \Gamma_D,\\
\boldsymbol{\sigma}(\boldsymbol{u})\boldsymbol{n} &= \bar{\boldsymbol{t}}
&&\text{on } \Gamma_N,
\end{aligned}
\end{equation}
where \(\boldsymbol{n}\) is the outward unit normal vector, and \(\partial\Omega=\Gamma_D\cup\Gamma_N\) with \(\Gamma_D\cap\Gamma_N=\emptyset\). Define the trial and test spaces
\begin{equation}
\begin{aligned}
\mathcal{V}_{u} &=
\left\{
\boldsymbol{u}\in H^1(\Omega;\mathbb{R}^3)
\,\middle|\,
\boldsymbol{u}=\bar{\boldsymbol{u}}
\text{ on }\Gamma_D
\right\},\\
\mathcal{V}_{0} &=
\left\{
\delta\boldsymbol{u}\in H^1(\Omega;\mathbb{R}^3)
\,\middle|\,
\delta\boldsymbol{u}=\boldsymbol{0}
\text{ on }\Gamma_D
\right\}.
\end{aligned}
\end{equation}

By multiplying the equilibrium equation by a test function \(\delta \boldsymbol{u}\), integrating over the domain, applying integration by parts, and using the boundary conditions, the weak form is obtained as
\begin{equation}
\text{find \(\boldsymbol{u}\in \mathcal{V}_{u}\) such that:}
\qquad
a(\boldsymbol{u},\delta \boldsymbol{u})=\ell( \delta\boldsymbol{u})
\qquad
\forall \delta \boldsymbol{u}\in \mathcal{V}_{0},
\end{equation}
where
\begin{equation}
a(\boldsymbol{u},\delta \boldsymbol{u})
=
\int_{\Omega}
\boldsymbol{\varepsilon}(\boldsymbol{u}):\mathbb{C}:\boldsymbol{\varepsilon}(\delta \boldsymbol{u})
\,d\Omega,
\end{equation}
and
\begin{equation}
\ell(\delta \boldsymbol{u}) = \int_{\Omega} \delta \boldsymbol{u}\cdot\boldsymbol{f}\,d\Omega +
\int_{\Gamma_N} \delta \boldsymbol{u}\cdot\bar{\boldsymbol{t}}\,d\Gamma.
\end{equation}
In this work, isogeometric discretization with NURBS/B-spline basis functions $\boldsymbol{N}$ is employed to approximate the displacement field as
\begin{equation}
\boldsymbol{u}_h = \boldsymbol{N}^T\boldsymbol{d},
\end{equation}
where $\boldsymbol{d}$ contains the discrete displacement coefficients associated with the spline basis. In a single-patch IGA domain, the multivariate shape function $\boldsymbol{N}=\boldsymbol{N}(\xi_1,\xi_2,\xi_3)$ is defined as $N_{ijk}(\xi_1,\xi_2,\xi_3) = N^1_i (\xi_1) N^2_j (\xi_2) N^3_k (\xi_3)$, where each $N^i(\xi_i)$ is a univariate shape function. The application of IGA to elasticity problems has been well established in the literature. Therefore, we omit the detailed construction of B-spline/NURBS basis functions and refinement techniques and refer interested readers to the standard references \cite{hughes2005isogeometric, piegl1997nurbs, cottrell2007studies, nguyen2015isogeometric}. Introducing the strain-displacement matrix $\boldsymbol{B}$, the discrete stiffness matrix and force vector are assembled from
\begin{equation}
\boldsymbol{K}
=
\int_{\Omega} \boldsymbol{B}^{T} \boldsymbol{C} \boldsymbol{B}\,d\Omega,
\end{equation}
and
\begin{equation}
\boldsymbol{f}
=
\int_{\Omega} \boldsymbol{N}^T \boldsymbol{f}\,d\Omega
+
\int_{\Gamma_N} \boldsymbol{N}^T \bar{\boldsymbol{t}}\,d\Gamma.
\end{equation}
Here, the differential operator $\boldsymbol{B}$ at a quadrature point $(i,j,k)$ is
\begin{equation}
    \boldsymbol{B}_{ijk} = \begin{bmatrix}
        \frac{\partial N_{ijk}}{\partial x} &0 &0\\
        0 &\frac{\partial N_{ijk}}{\partial y} &0\\
        0 &0 &\frac{\partial N_{ijk}}{\partial z}\\
        0 &\frac{\partial N_{ijk}}{\partial z} &\frac{\partial N_{ijk}}{\partial y}\\
        \frac{\partial N_{ijk}}{\partial z} &0 &\frac{\partial N_{ijk}}{\partial x}\\        
        \frac{\partial N_{ijk}}{\partial y} &\frac{\partial N_{ijk}}{\partial x} &0
    \end{bmatrix}.
\end{equation}
For isotropic elasticity in 3D, the material tensor is
\begin{equation}
\boldsymbol{C}=
\begin{bmatrix}
\lambda+2\mu & \lambda & \lambda & 0 & 0 & 0\\
\lambda & \lambda+2\mu & \lambda & 0 & 0 & 0\\
\lambda & \lambda & \lambda+2\mu & 0 & 0 & 0\\
0 & 0 & 0 & \mu & 0 & 0\\
0 & 0 & 0 & 0 & \mu & 0\\
0 & 0 & 0 & 0 & 0 & \mu
\end{bmatrix}.
\end{equation}
Here, Voigt notation is used to simplify the derivation of the geometry-derived coefficient expressions. An equivalent formulation written in full tensor notation would lead to the same results. The following sections rewrite the standard elasticity operators in the reference coordinate system and then exploit their tensor-product structure for tensor-train assembly and solution.
\section{Tensor-train formulation for 3D elasticity}
To build a TT framework for elasticity problems, we express the weak-form operators as tensor-product expressions. To do so, we first define the matrix
\begin{equation}
\boldsymbol{S} =
\begin{bmatrix}
1 & 0 & 0 & 0 & 0 & 0 & 0 & 0 & 0\\
0 & 0 & 0 & 0 & 1 & 0 & 0 & 0 & 0\\
0 & 0 & 0 & 0 & 0 & 0 & 0 & 0 & 1\\
0 & 0 & 0 & 0 & 0 & 1 & 0 & 1 & 0\\
0 & 0 & 1 & 0 & 0 & 0 & 1 & 0 & 0\\
0 & 1 & 0 & 1 & 0 & 0 & 0 & 0 & 0
\end{bmatrix}.
\end{equation}
The reference derivative operator is defined as
\begin{equation}
\boldsymbol{D}_{ijk}=
\begin{bmatrix}
\dfrac{\partial N_{ijk}}{\partial \xi_1} &0 &0\\
\dfrac{\partial N_{ijk}}{\partial \xi_2} &0 &0\\
\dfrac{\partial N_{ijk}}{\partial \xi_3} &0 &0\\
0 &\dfrac{\partial N_{ijk}}{\partial \xi_1} &0\\
0 &\dfrac{\partial N_{ijk}}{\partial \xi_2} &0\\
0 &\dfrac{\partial N_{ijk}}{\partial \xi_3} &0\\
0 &0 &\dfrac{\partial N_{ijk}}{\partial \xi_1}\\
0 &0 &\dfrac{\partial N_{ijk}}{\partial \xi_2}\\
0 &0 &\dfrac{\partial N_{ijk}}{\partial \xi_3}
\end{bmatrix},
\end{equation}
The Jacobian matrix is defined as the derivative of the physical coordinates $\boldsymbol{x} =\{x_1, x_2, x_3\}$ with respect to the parametric coordinates $\boldsymbol{\xi} = \{\xi_1, \xi_2, \xi_3\}$ as
\begin{equation}
    \boldsymbol{J} = \begin{bmatrix}
        \frac{\partial x_1}{\partial \xi_1} &\frac{\partial x_1}{\partial \xi_2} &\frac{\partial x_1}{\partial \xi_3}\\
        \frac{\partial x_2}{\partial \xi_1} &\frac{\partial x_2}{\partial \xi_2} &\frac{\partial x_2}{\partial \xi_3}\\ \frac{\partial x_3}{\partial \xi_1} &\frac{\partial x_3}{\partial \xi_2} &\frac{\partial x_3}{\partial \xi_3}
    \end{bmatrix}, \qquad \text{its inverse }
    \boldsymbol{J}^{-1} = \begin{bmatrix}
        \frac{\partial \xi_1}{\partial x_1} &\frac{\partial \xi_1}{\partial x_2} &\frac{\partial \xi_1}{\partial x_3}\\
        \frac{\partial \xi_2}{\partial x_1} &\frac{\partial \xi_2}{\partial x_2} &\frac{\partial \xi_2}{\partial x_3}\\ \frac{\partial \xi_3}{\partial x_1} &\frac{\partial \xi_3}{\partial x_2} &\frac{\partial \xi_3}{\partial x_3}
        \end{bmatrix},\qquad\text{ and}
\end{equation}
and the Jacobian determinant is $|J| = \det(\boldsymbol{J})$.
The expanded Jacobian matrix and its inverse are expressed as
\begin{equation}
    \boldsymbol{J}_{e} = \begin{bmatrix}
        \boldsymbol{J} &\boldsymbol{0} &\boldsymbol{0}\\
        \boldsymbol{0} &\boldsymbol{J} &\boldsymbol{0}\\
        \boldsymbol{0} &\boldsymbol{0} &\boldsymbol{J}
    \end{bmatrix},
    \qquad\text{ and }\qquad 
        \boldsymbol{J}^{-1}_{e} = \begin{bmatrix}
        \boldsymbol{J}^{-1} &\boldsymbol{0} &\boldsymbol{0}\\
        \boldsymbol{0} &\boldsymbol{J}^{-1} &\boldsymbol{0}\\
        \boldsymbol{0} &\boldsymbol{0} &\boldsymbol{J}^{-1}
    \end{bmatrix},
\end{equation}
The differential operator $\boldsymbol{B}$ is expressed as
\begin{equation}
\boldsymbol{B} = \boldsymbol{S} (\boldsymbol{J}_{e}^{-T} \boldsymbol{D}).
\end{equation}
The Voigt notation is encoded through the matrix $\boldsymbol{S}$; this allows us to define
\begin{equation}
\boldsymbol{R}
=
\boldsymbol{J}^{-1}_{e}\, \boldsymbol{S}^T\, \boldsymbol{C} \boldsymbol{S} \boldsymbol{J}_{e}^{-T}\, |J|.
\end{equation}
At a quadrature point $(i,j,k)$ in 3D, $\boldsymbol{R}$ is a $9 \times 9$ matrix
\begin{equation}
\boldsymbol{R} =
\begin{bmatrix}
R_{11}&R_{12}&R_{13}&R_{14}&R_{15}&R_{16}&R_{17}&R_{18}&R_{19}\\
R_{21}&R_{22}&R_{23}&R_{24}&R_{25}&R_{26}&R_{27}&R_{28}&R_{29}\\
R_{31}&R_{32}&R_{33}&R_{34}&R_{35}&R_{36}&R_{37}&R_{38}&R_{39}\\
R_{41}&R_{42}&R_{43}&R_{44}&R_{45}&R_{46}&R_{47}&R_{48}&R_{49}\\
R_{51}&R_{52}&R_{53}&R_{54}&R_{55}&R_{56}&R_{57}&R_{58}&R_{59}\\
R_{61}&R_{62}&R_{63}&R_{64}&R_{65}&R_{66}&R_{67}&R_{68}&R_{69}\\
R_{71}&R_{72}&R_{73}&R_{74}&R_{75}&R_{76}&R_{77}&R_{78}&R_{79}\\
R_{81}&R_{82}&R_{83}&R_{84}&R_{85}&R_{86}&R_{87}&R_{88}&R_{89}\\
R_{91}&R_{92}&R_{93}&R_{94}&R_{95}&R_{96}&R_{97}&R_{98}&R_{99}
\end{bmatrix}.
\end{equation}
A few entries are
\begin{equation}
\begin{split}
R_{11} = R_{11}(\xi_1,\xi_2,\xi_3) &= |J|\left[
\left(\frac{\partial \xi_1}{\partial x_2}\right)^2 \mu
+
\left(\frac{\partial \xi_1}{\partial x_3}\right)^2 \mu
+
\left(\frac{\partial \xi_1}{\partial x_1}\right)^2(\lambda+2\mu)
\right],\\
R_{12} = R_{12}(\xi_1,\xi_2,\xi_3) &= |J|\left[
\frac{\partial \xi_1}{\partial x_2}\frac{\partial \xi_2}{\partial x_2}\mu
+
\frac{\partial \xi_1}{\partial x_3}\frac{\partial \xi_2}{\partial x_3}\mu
+
\frac{\partial \xi_1}{\partial x_1}\frac{\partial \xi_2}{\partial x_1}(\lambda+2\mu)
\right].
\end{split}    
\end{equation}
Each component of $\boldsymbol{R}$ contains geometry-derived coefficient fields (metric tensor and Jacobian determinant). The stiffness matrix is expressed as
\begin{equation}
\boldsymbol{K} = \int_0^1\int_0^1\int_0^1 \boldsymbol{D}^T \boldsymbol{R} \boldsymbol{D} d\xi_1d\xi_2d\xi_3,
\end{equation}
To expose the displacement-component structure used in the tensor-product assembly, it is useful to partition the \(9\times 9\) matrix \(\boldsymbol{R}\) into nine \(3\times 3\) component submatrices:
\begin{equation}
\boldsymbol{R} =
\begin{bmatrix}
\boldsymbol{R}_{xx} &\boldsymbol{R}_{xy} &\boldsymbol{R}_{xz}\\
\boldsymbol{R}_{yx} &\boldsymbol{R}_{yy} &\boldsymbol{R}_{yz}\\
\boldsymbol{R}_{zx} &\boldsymbol{R}_{zy} &\boldsymbol{R}_{zz}
\end{bmatrix},
\end{equation}
where
\begin{equation}
\boldsymbol{R}_{xx}=
\begin{bmatrix}
R_{11}&R_{12}&R_{13}\\
R_{21}&R_{22}&R_{23}\\
R_{31}&R_{32}&R_{33}
\end{bmatrix},
\qquad
\boldsymbol{R}_{xy}=
\begin{bmatrix}
R_{14}&R_{15}&R_{16}\\
R_{24}&R_{25}&R_{26}\\
R_{34}&R_{35}&R_{36}
\end{bmatrix},
\qquad
\boldsymbol{R}_{xz}=
\begin{bmatrix}
R_{17}&R_{18}&R_{19}\\
R_{27}&R_{28}&R_{29}\\
R_{37}&R_{38}&R_{39}
\end{bmatrix},
\end{equation}
\begin{equation}
\boldsymbol{R}_{yx}=
\begin{bmatrix}
R_{41}&R_{42}&R_{43}\\
R_{51}&R_{52}&R_{53}\\
R_{61}&R_{62}&R_{63}
\end{bmatrix},
\qquad
\boldsymbol{R}_{yy}=
\begin{bmatrix}
R_{44}&R_{45}&R_{46}\\
R_{54}&R_{55}&R_{56}\\
R_{64}&R_{65}&R_{66}
\end{bmatrix},
\qquad
\boldsymbol{R}_{yz}=
\begin{bmatrix}
R_{47}&R_{48}&R_{49}\\
R_{57}&R_{58}&R_{59}\\
R_{67}&R_{68}&R_{69}
\end{bmatrix},
\end{equation}
\begin{equation}
\boldsymbol{R}_{zx}=
\begin{bmatrix}
R_{71}&R_{72}&R_{73}\\
R_{81}&R_{82}&R_{83}\\
R_{91}&R_{92}&R_{93}
\end{bmatrix},
\qquad
\boldsymbol{R}_{zy}=
\begin{bmatrix}
R_{74}&R_{75}&R_{76}\\
R_{84}&R_{85}&R_{86}\\
R_{94}&R_{95}&R_{96}
\end{bmatrix},
\qquad
\boldsymbol{R}_{zz}=
\begin{bmatrix}
R_{77}&R_{78}&R_{79}\\
R_{87}&R_{88}&R_{89}\\
R_{97}&R_{98}&R_{99}
\end{bmatrix}.
\end{equation}
Let the one-dimensional matrix factors be denoted by
\begin{equation}
\begin{split}
    \boldsymbol{\phi}^{(1)}_{xx} = \frac{\partial \boldsymbol{N}_1}{\partial \xi_1} \frac{\partial \boldsymbol{N}_1^T}{\partial \xi_1},\qquad
    \boldsymbol{\phi}^{(2)}_{xx} = \boldsymbol{N}_2\boldsymbol{N}_2^T,\qquad
    \boldsymbol{\phi}^{(3)}_{xx} = \boldsymbol{N}_3\boldsymbol{N}_3^T,\\
    \boldsymbol{\phi}^{(1)}_{xy} = \frac{\partial \boldsymbol{N}_1}{\partial \xi_1}\boldsymbol{N}_1^T,\qquad
    \boldsymbol{\phi}^{(2)}_{xy} = \boldsymbol{N}_2\frac{\partial \boldsymbol{N}_2^T}{\partial \xi_2},\qquad
    \boldsymbol{\phi}^{(3)}_{xy} = \boldsymbol{N}_3\boldsymbol{N}_3^T,\\
    \boldsymbol{\phi}^{(1)}_{xz} = \frac{\partial \boldsymbol{N}_1}{\partial \xi_1}\boldsymbol{N}_1^T,\qquad
    \boldsymbol{\phi}^{(2)}_{xz} = \boldsymbol{N}_2\boldsymbol{N}_2^T,\qquad
    \boldsymbol{\phi}^{(3)}_{xz} = \boldsymbol{N}_3\frac{\partial \boldsymbol{N}_3^T}{\partial \xi_3},\\
    \boldsymbol{\phi}^{(1)}_{yx} = \boldsymbol{N}_1\frac{\partial \boldsymbol{N}_1^T}{\partial \xi_1},\qquad
    \boldsymbol{\phi}^{(2)}_{yx} = \frac{\partial \boldsymbol{N}_2}{\partial \xi_2}\boldsymbol{N}_2^T,\qquad
    \boldsymbol{\phi}^{(3)}_{yx} = \boldsymbol{N}_3\boldsymbol{N}_3^T,\\
    \boldsymbol{\phi}^{(1)}_{yy} = \boldsymbol{N}_1\boldsymbol{N}_1^T,\qquad
    \boldsymbol{\phi}^{(2)}_{yy} = \frac{\partial \boldsymbol{N}_2}{\partial \xi_2}\frac{\partial \boldsymbol{N}_2^T}{\partial \xi_2},\qquad
    \boldsymbol{\phi}^{(3)}_{yy} = \boldsymbol{N}_3\boldsymbol{N}_3^T,\\
    \boldsymbol{\phi}^{(1)}_{yz} = \boldsymbol{N}_1\boldsymbol{N}_1^T,\qquad
    \boldsymbol{\phi}^{(2)}_{yz} = \frac{\partial \boldsymbol{N}_2}{\partial \xi_2}\boldsymbol{N}_2^T,\qquad
    \boldsymbol{\phi}^{(3)}_{yz} = \boldsymbol{N}_3\frac{\partial \boldsymbol{N}_3^T}{\partial \xi_3},\\
    \boldsymbol{\phi}^{(1)}_{zx} = \boldsymbol{N}_1\frac{\partial \boldsymbol{N}_1^T}{\partial \xi_1},\qquad
    \boldsymbol{\phi}^{(2)}_{zx} = \boldsymbol{N}_2\boldsymbol{N}_2^T,\qquad
    \boldsymbol{\phi}^{(3)}_{zx} = \frac{\partial \boldsymbol{N}_3}{\partial \xi_3}\boldsymbol{N}_3^T,\\
    \boldsymbol{\phi}^{(1)}_{zy} = \boldsymbol{N}_1\boldsymbol{N}_1^T,\qquad
    \boldsymbol{\phi}^{(2)}_{zy} = \boldsymbol{N}_2\frac{\partial \boldsymbol{N}_2^T}{\partial \xi_2},\qquad
    \boldsymbol{\phi}^{(3)}_{zy} = \frac{\partial \boldsymbol{N}_3}{\partial \xi_3}\boldsymbol{N}_3^T,\\
    \boldsymbol{\phi}^{(1)}_{zz} = \boldsymbol{N}_1\boldsymbol{N}_1^T,\qquad
    \boldsymbol{\phi}^{(2)}_{zz} = \boldsymbol{N}_2\boldsymbol{N}_2^T,\qquad
    \boldsymbol{\phi}^{(3)}_{zz} =\frac{\partial \boldsymbol{N}_3}{\partial \xi_3}\frac{\partial \boldsymbol{N}_3^T}{\partial \xi_3}.
\end{split}
\label{eq.phi_factors}
\end{equation}
The displacement-component stiffness operators can be written as
\begin{equation}
\boldsymbol{K}_{ab}
=
\int_0^1
\int_0^1
\int_0^1
\sum_{c=1}^{3}
\sum_{d=1}^{3}
R_{abcd}
\left(
\boldsymbol{\phi}^{(1)}_{cd}
\otimes
\boldsymbol{\phi}^{(2)}_{cd}
\otimes
\boldsymbol{\phi}^{(3)}_{cd}
\right)
\,d\xi_1d\xi_2d\xi_3,
\label{eq.K9_compact}
\end{equation}
Each component $\boldsymbol{K}_{ab} \in \mathbb{R}^{(n_1n_2n_3) \times (n_1n_2n_3)}$, where \(a,b\in\{x,y,z\}\) denote the displacement-component indices and \(c,d\in\{1,2,3\}\) denote the local row and column indices inside the corresponding \(3\times 3\) coefficient submatrix
\(\boldsymbol{R}_{ab}\). The matrix $\boldsymbol{K} \in \mathbb{R}^{(3n_1n_2n_3) \times (3n_1n_2n_3)} $ is written in displacement-component form as
\begin{equation}
\boldsymbol{K} = 
\begin{bmatrix}
\boldsymbol{K}_{xx} & \boldsymbol{K}_{xy} & \boldsymbol{K}_{xz}\\
\boldsymbol{K}_{yx} & \boldsymbol{K}_{yy} & \boldsymbol{K}_{yz}\\
\boldsymbol{K}_{zx} & \boldsymbol{K}_{zy} & \boldsymbol{K}_{zz}
\end{bmatrix},
\end{equation}
which is assembled as
\begin{equation}
\boldsymbol{K}
=
\sum_{a,b \in \{x,y,z\}}
\boldsymbol{K}_{ab}\otimes \boldsymbol{E}_{ab},
\label{eq.K_component_assembly}
\end{equation}
where $\boldsymbol{E}_{ab}\in\mathbb{R}^{3\times 3}$ is the component-selection matrix with one nonzero entry at position $(a,b)$. 
The tensor-train representation of $\boldsymbol{K}$ is obtained after tensorizing the row and column indices according to the prescribed TT ordering. The volumetric force is
\begin{equation}
    \boldsymbol{f} = \begin{bmatrix}
        \boldsymbol{f}_x\\
        \boldsymbol{f}_y\\
        \boldsymbol{f}_z
    \end{bmatrix},
\label{eq.fdec}
\end{equation}
where
\begin{equation}
\begin{split}
    &\boldsymbol{f}_x = \int_0^1\int_0^1 \int_0^1f_x(x,y,z)\boldsymbol{N}_1^T\otimes\boldsymbol{N}_2^T\otimes\boldsymbol{N}_3^T |J|d\xi_1d\xi_2d\xi_3,\\
    &\boldsymbol{f}_y = \int_0^1\int_0^1 \int_0^1f_y(x,y,z)\boldsymbol{N}_1^T\otimes\boldsymbol{N}_2^T\otimes\boldsymbol{N}_3^T |J|d\xi_1d\xi_2d\xi_3,\\
    &\boldsymbol{f}_z = \int_0^1\int_0^1 \int_0^1f_z(x,y,z)\boldsymbol{N}_1^T\otimes\boldsymbol{N}_2^T\otimes\boldsymbol{N}_3^T |J| d\xi_1d\xi_2d\xi_3.\\
\end{split}
\end{equation}
\begin{figure}
\begin{center}
\includegraphics[width=0.99\textwidth]{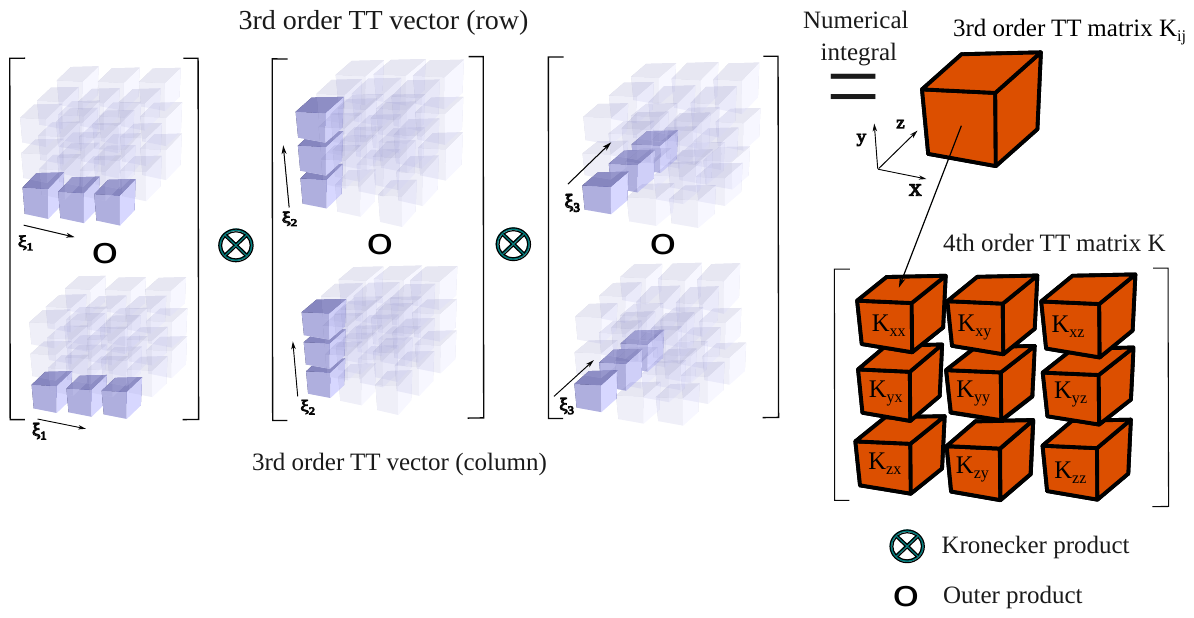}
\end{center}
\caption{Construction of the coupled fourth-order elasticity TT matrix $\boldsymbol{K}^{\mathtt{TT}}$. The geometry-dependent coefficient fields $R_{abcd}$ are first represented in the three spatial TT modes and used to construct the corresponding displacement-component operators $\boldsymbol{K}_{ab}^{\mathtt{TT}}$. A fourth TT core then introduces the displacement-component indices and forms the complete coupled operator.}
\label{Fig.Constructing_K4}
\end{figure}
\section{TT construction of the elasticity operator}
The geometry-dependent coefficient $R_{abcd}(\boldsymbol{\xi})$ and the Jacobian determinant $|J(\boldsymbol{\xi})|$ are spatially varying scalars over the parametric domain. To construct the TT representation, the Jacobian determinant $|J(\boldsymbol{\xi})|$ is sampled at 3D quadrature points \((i_1,i_2,i_3)\) using TT-cross interpolation \cite{oseledets2010tt} with a prescribed accuracy \(\varepsilon\) to obtain the three-dimensional array $\left(|J|\right)_{i_1 i_2 i_3}$ as
\begin{equation}
\left(|J|\right)_{i_1 i_2 i_3}
\approx
\left(|J|^{\mathtt{TT}}\right)_{i_1 i_2 i_3}
=
Q^{(1)}_{1i_1\alpha_1}
Q^{(2)}_{\alpha_1 i_2 \alpha_2}
Q^{(3)}_{\alpha_2 i_3 1},
\label{eq.Jtt}
\end{equation}
where $Q^{(i)}$ are the cores associated with the $\xi_i$ direction, and the repeated TT-rank indices \(\alpha_1\) and \(\alpha_2\) imply Einstein summation. The TT representation of \(|J|\) is used in the construction of the volumetric force vector. The Neumann surface force is constructed separately using the corresponding surface Jacobian, as described in \Aref{App.TTNeumann}.

Similarly, the geometry-dependent coefficient fields are also sampled to obtain $R_{abcd}(i_1i_2i_3)$, or equivalently $\left(R_{abcd}\right)_{i_1i_2i_3}$, as
\begin{equation}
\left(R_{a b c d}\right)_{i_1 i_2 i_3}
\approx
\left(R^{\mathtt{TT}}_{a b c d}\right)_{i_1 i_2 i_3}
=
G^{abcd(1)}_{1 i_1 \alpha_1}
G^{abcd(2)}_{\alpha_1 i_2 \alpha_2}
G^{abcd(3)}_{\alpha_2 i_3 1}.
\label{eq:Rtt}
\end{equation}
Here, $G^{abcd(i)}$ are the TT cores of $R_{abcd}$, where $(a,b\in{x,y,z})$ denote the displacement-component indices and $(c,d\in{1,2,3})$ denote the derivative-direction indices. These component indices are not summed in \eqref{eq:Rtt}. The TT approximation of $R_{abcd}$ from \Eqref{eq:Rtt} is used to construct the TT representation of each displacement-component stiffness operator \(\boldsymbol{K}_{ab}\). It should be emphasized that the TT cores of \(R_{abcd}\) are not directly the TT cores of \(\boldsymbol{K}_{ab}\). Instead, each one-dimensional core of \(R_{abcd}^{\mathtt{TT}}\) provides a set of weights used to assemble the corresponding one-dimensional stiffness matrices. For a fixed coefficient \(R_{abcd}^{\mathtt{TT}}\), the TT cores in \Eqref{eq:Rtt} are used as one-dimensional quadrature weights to compute the TT core of the stiffness matrix
\begin{equation}
\begin{split}
&T^{abcd(1)}_{1 p_1 q_1 \alpha_1}
= \int_0^1G^{abcd(1)}_{1 i_1 \alpha_1}\Big(\boldsymbol{\phi}^{(1)}_{cd}\Big)_{p_1q_1}\,d\xi_1 
= \sum_{s_1=1}^{n_{q_1}}w^{(1)}_{s_1}G^{abcd(1)}_{1 s_1 \alpha_1}
\Big(\boldsymbol{\phi}^{(1)}_{cd}(\xi_{s_1})\Big)_{p_1q_1},
\\
&T^{abcd(2)}_{\alpha_1p_2q_2\alpha_2}
=\int_0^1 G^{abcd(2)}_{\alpha_1 i_2 \alpha_2}\Big(\boldsymbol{\phi}^{(2)}_{cd}\Big)_{p_2q_2}\,d\xi_2 
=\sum_{s_2=1}^{n_{q_2}}w^{(2)}_{s_2}G^{abcd(2)}_{\alpha_1 s_2 \alpha_2}\Big(\boldsymbol{\phi}^{(2)}_{cd}(\xi_{s_2})\Big)_{p_2q_2},
\\
&T^{abcd(3)}_{\alpha_2 p_3q_3 1}
=\int_0^1 G^{abcd(3)}_{\alpha_2 i_3 1}\Big(\boldsymbol{\phi}^{(3)}_{cd}\Big)_{p_3q_3}d\xi_3
=\sum_{s_3=1}^{n_{q_3}}w^{(3)}_{s_3}G^{abcd(3)}_{\alpha_2 s_3 1}\Big(\boldsymbol{\phi}^{(3)}_{cd}(\xi_{s_3})\Big)_{p_3q_3},
\end{split}
\label{eq:T-from-R-cores}
\end{equation}
where \(\xi_{s_\ell}\) and \(w_{s_\ell}\) denote the quadrature point and the corresponding quadrature weight in the \(\xi_\ell\)-direction, respectively, and \(n_{q_\ell}\) is the number of quadrature points used in that direction. Using these weighted one-dimensional matrices, the TT representation of each displacement-component stiffness operator is constructed in index form as
\begin{equation}
\Big(K_{ab}\Big)_{p_1q_1p_2q_2p_3q_3} \approx \Big(K^{\mathtt{TT}}_{ab}\Big)_{p_1q_1p_2q_2p_3q_3}
=
\sum_{c=1}^{3}
\sum_{d=1}^{3}
T^{abcd(1)}_{1p_1q_1\alpha_1}
T^{abcd(2)}_{\alpha_1p_2q_2\alpha_2}
T^{abcd(3)}_{\alpha_2p_3q_31},
\label{eq:Kab0-TT-index}
\end{equation}
or it can be expressed as
\begin{equation}
\boldsymbol{K}^{\mathtt{TT}}_{ab}
=
\sum_{c=1}^{3}
\sum_{d=1}^{3}
\boldsymbol{T}^{abcd(1)}\otimes
\boldsymbol{T}^{abcd(2)}\otimes
\boldsymbol{T}^{abcd(3)}.
\label{eq:Kab-TT-index}
\end{equation}
Each term in \Eqref{eq:Kab-TT-index} corresponds to a Kronecker product of three weighted one-dimensional matrices. The summation over \(c\) and \(d\) accumulates the nine derivative-direction contributions associated with \(\boldsymbol{K}_{ab}\). The symmetry of the elasticity operator is also preserved in the constructed TT representation. Therefore, only the six upper-triangular displacement-component operators need to be constructed and stored through their TT cores; the remaining three are obtained by transposition:
\begin{equation}
\boldsymbol{K}_{yx}^{\mathtt{TT}}=(\boldsymbol{K}_{xy}^{\mathtt{TT}})^T,\qquad
\boldsymbol{K}_{zx}^{\mathtt{TT}}=(\boldsymbol{K}_{xz}^{\mathtt{TT}})^T,\qquad
\boldsymbol{K}_{zy}^{\mathtt{TT}}=(\boldsymbol{K}_{yz}^{\mathtt{TT}})^T.    
\end{equation}
The complete coupled stiffness matrix is represented as a fourth-order TT-matrix,
\begin{equation}
\big(K\big)_{p_1 q_1 p_2 q_2 p_3 q_3 p_4 q_4}
\approx
\big(K^{\mathtt{TT}}\big)_{p_1 q_1 p_2 q_2 p_3 q_3 p_4 q_4}
=
\sum_{a,b \in \{x,y,z\}} \Big(K^{\mathtt{TT}}_{ab}\Big)_{p_1 q_1 p_2 q_2 p_3 q_3}G^{ab(4)}_{1 p_4 q_4 1}
\label{eq.Ktt_4th}
\end{equation}
where \(p_4,q_4 \in \{1,2,3\}\) denote the displacement-component indices, and the fourth core is defined by
\begin{equation}
G^{ab(4)}_{1 p_4 q_4 1} = (\boldsymbol{E}_{ab})_{p_4 q_4}.
\end{equation}
The TT representation given in \eqref{eq.Ktt_4th} separates the spatial directions and the displacement-component direction into different TT cores. As a result, the information associated with each direction is stored in a single core. The construction process for \(\boldsymbol{K}^{\mathtt{TT}}\) is illustrated in \Figref{Fig.Constructing_K4}. Imposing boundary conditions in TT format requires modifying the corresponding TT core; this procedure is discussed in \Aref{Sec.BCs}. 

\section{Solution of the TT elasticity system with AMEn}
\label{Sec.SolvingTT}

After construction of the coupled stiffness operator and application of the boundary conditions, the discrete elasticity problem is represented entirely in TT format as
\begin{equation}
\boldsymbol{K}^{\mathtt{TT}}
\boldsymbol{u}^{\mathtt{TT}}
=
\boldsymbol{f}^{\mathtt{TT}}.
\end{equation}
The displacement-component index is included as an additional TT mode, so the three displacement components are represented and solved simultaneously in one coupled TT system. The system is solved using the alternating minimal energy (AMEn) method \cite{dolgov2014alternating}. AMEn performs alternating optimization over the TT cores of the displacement field and enriches the approximation using residual information. The resulting ranks are therefore adapted during the solution process.

Three numerical tolerances are distinguished in the implementation. The TT-cross interpolation tolerance, \(\varepsilon_{\mathrm{cross}}\), controls the approximation of the sampled geometry-dependent coefficient fields; the TT truncation tolerance, \(\varepsilon_{\mathrm{TT}}\), controls rank reduction during TT algebraic operations; and the AMEn tolerance, \(\varepsilon_{\mathrm{AMEn}}\), controls convergence through the default residual-based criterion, which monitors the maximum relative residual of the local projected problems over an AMEn sweep. In the numerical examples, all three prescribed tolerances are set to \(10^{-8}\). For comparison, the conventional full-grid IGA formulation explicitly assembles the sparse stiffness matrix and solves the reduced linear system using a sparse direct solver. The full-grid matrix storage is measured by the number of nonzero entries, whereas TT storage is measured by the total number of scalar entries contained in the TT cores.
All computations were performed in MATLAB R2022b using TT-Toolbox 2.2 on a Mac workstation equipped with an Apple M1 Max processor and 32~GB of unified memory. All calculations were performed in serial. For the conventional full-grid IGA formulation, the reported computational time includes sparse-matrix assembly and solution of the resulting linear system using the sparse direct solver. For the TT formulation, the reported computational time includes TT-cross approximation of the geometry-dependent coefficient fields, construction of the TT operator and force vector, and solution of the resulting TT system using AMEn. Thus, the reported timings represent the complete assembly-and-solution workflow for both formulations.
\section{Numerical examples}
The material parameters are chosen as \(\lambda=10^{7}\,\mathrm{MPa}\) and \(\mu=0.7\times 10^{7}\,\mathrm{MPa}\). Storage efficiency is measured using the actual number of stored scalar entries. For the assembled full-grid stiffness matrix, this quantity is its number of nonzero entries (NNZ); for TT matrices and vectors, it is the total number of scalar entries contained in the TT cores. The corresponding storage and compression-ratio definitions are given in \Aref{sec.Appen-Storage}. The cylindrical example provides a direct verification against the full-grid discretization in terms of accuracy, storage, and solution time. The open-hemisphere and square-hole examples then assess compression and computational scalability over substantially larger problem sizes, while full-grid timing results are reported only over the range for which explicit sparse assembly and solution remain feasible.

\subsection{Cylindrical structure}
\label{sec.Cylindrical}
\begin{figure}
\begin{center}
\includegraphics[width=0.60\textwidth]{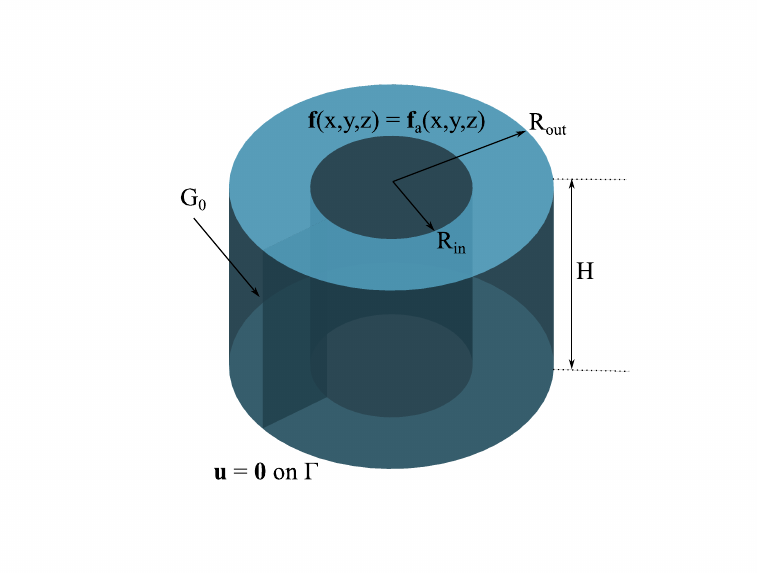}
\end{center}
\caption{Cylindrical structure subjected to the boundary condition $\boldsymbol{u} = \boldsymbol{0}$ on its surrounding surfaces $\Gamma$ and the volumetric force $f_a(x,y,z)$. An open knot vector is employed to generate the geometry, and a $G_0$ constraint is imposed for geometric continuity in the domain.}
\label{Fig.Ring}
\end{figure}
In this numerical example, a cylindrical structure subjected to a volumetric force is studied, as illustrated in \Figref{Fig.Ring}. The geometry is defined by \(R_{\mathrm{out}}=H=1\,\mathrm{m}\) and \(R_{\mathrm{in}}=0.5\,\mathrm{m}\). Quadratic basis functions are used in all three parametric directions, i.e., \(p_1=p_2=p_3=2\). A manufactured solution is used to verify the proposed formulation:
\begin{equation}
\boldsymbol{u}(x,y,z) = \phi(x,y,z)
\begin{bmatrix}
1 \ 2 \ 3
\end{bmatrix},
\end{equation}
where \(\phi(x,y,z)=P(x,y,z)S(x,y,z)\). The detailed expressions for \(P(x,y,z)\) and \(S(x,y,z)\) are provided in \Aref{sec.Appen-Manufacture-Ring}. The manufactured solution is constructed to satisfy the prescribed boundary conditions $\boldsymbol{u} = \boldsymbol{0}$ on $\Gamma$ and to yield vanishing natural traction. The corresponding volumetric force is then obtained by substituting this solution into the strong form of the governing equations. The detailed derivation of the volumetric force is also provided in \Aref{sec.Appen-Manufacture-Ring}.

\begin{figure}
\begin{center}
\minipage{0.5\textwidth}
\includegraphics[width=0.80\textwidth]{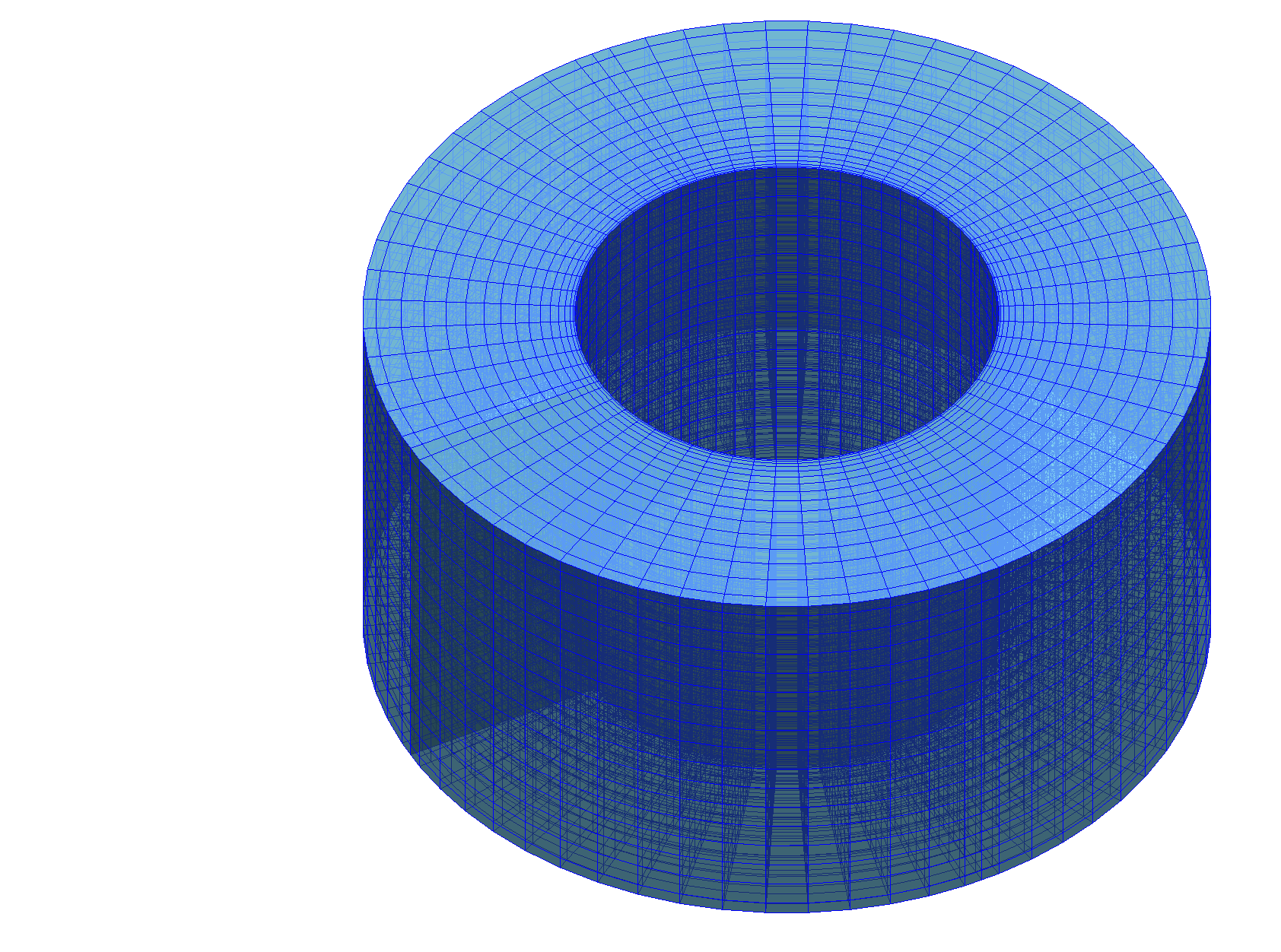}
\caption*{(a) Control mesh grid $69 \times 18 \times 18$}
\endminipage
\hfill
\minipage{0.5\textwidth}
\includegraphics[width=1.0\textwidth]{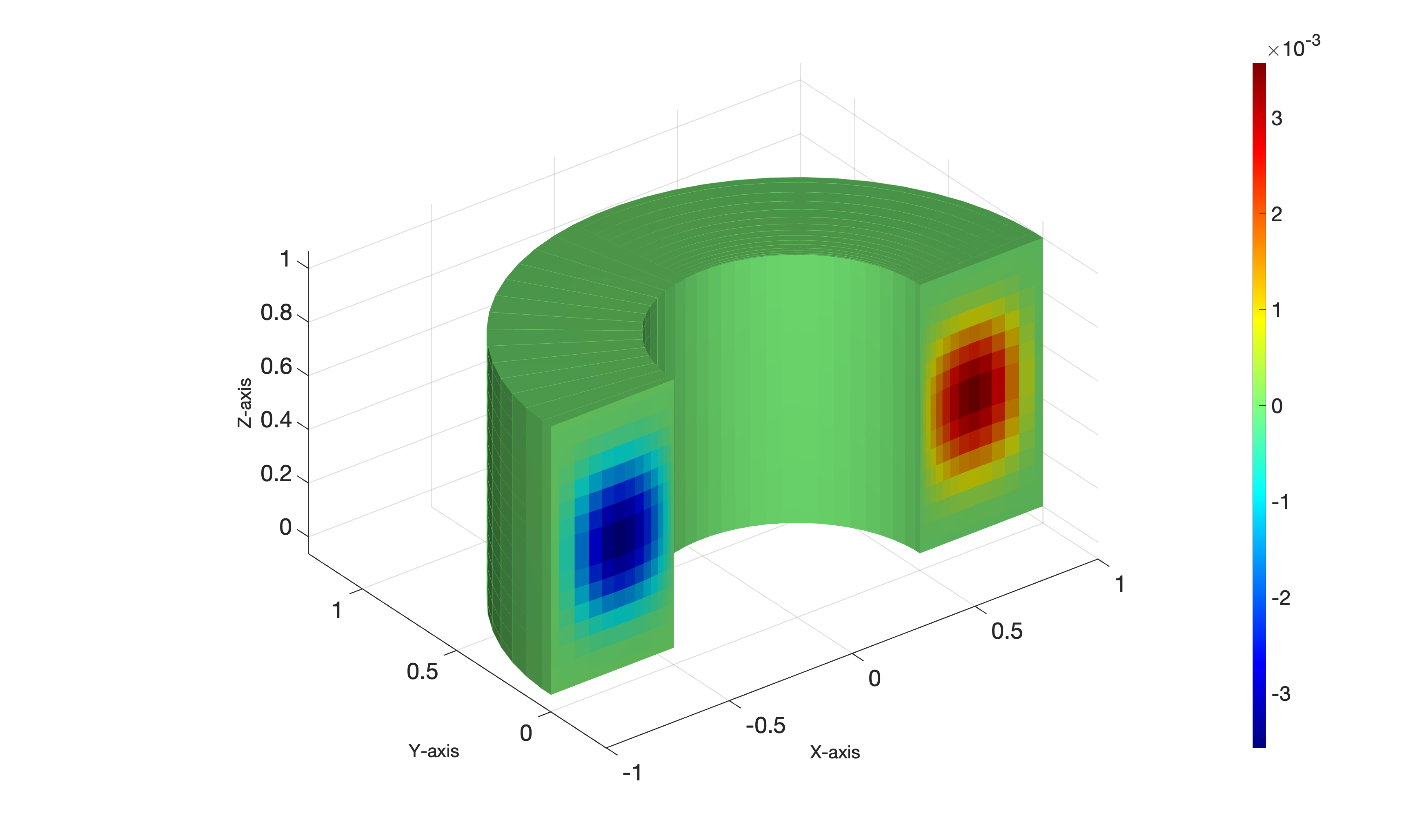}
\caption*{(b) $x$-displacement}
\endminipage
\vfill
\minipage{0.5\textwidth}
\includegraphics[width=1.0\textwidth]{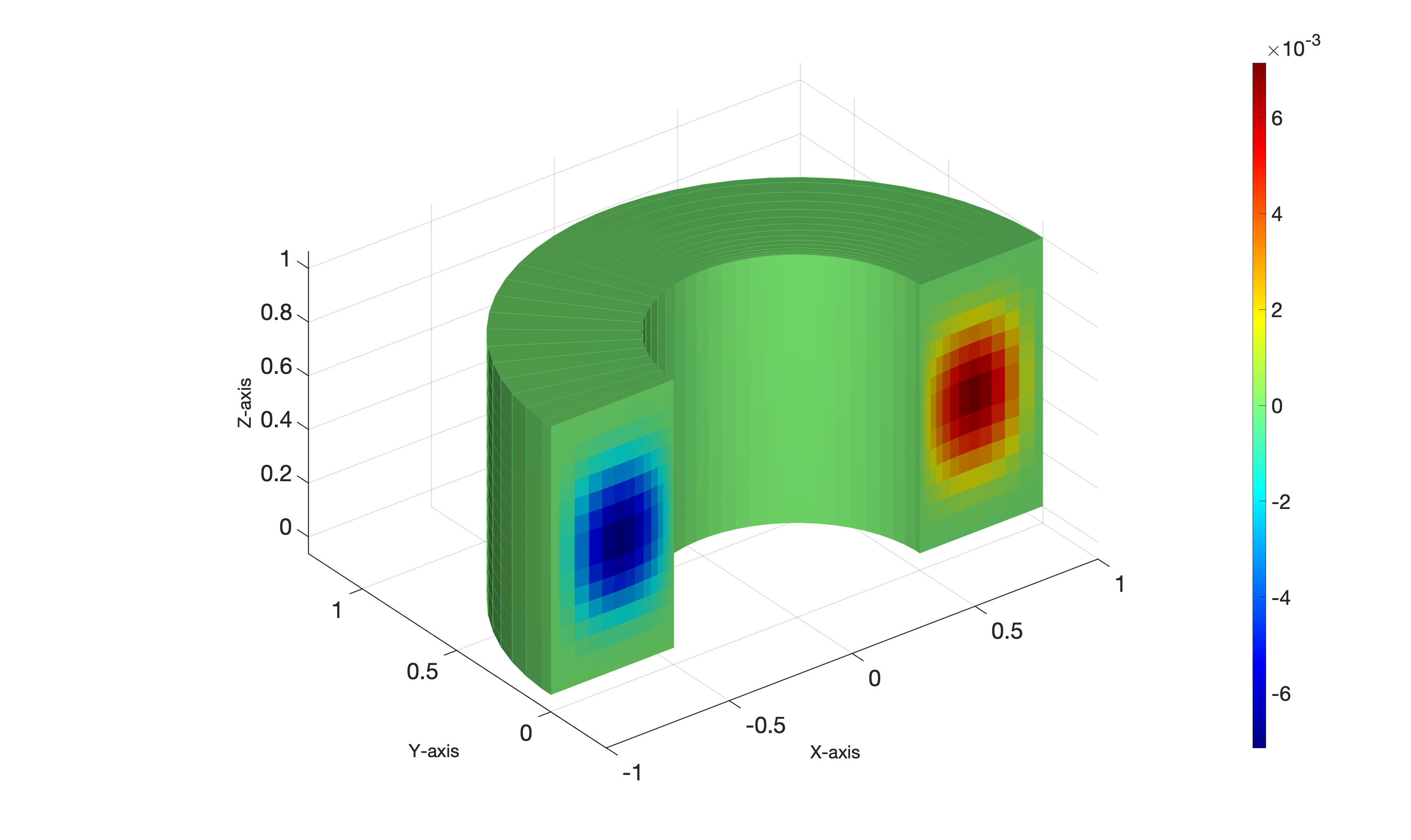}
\caption*{(c) $y$-displacement}
\endminipage
\hfill
\minipage{0.5\textwidth}
\includegraphics[width=1.0\textwidth]{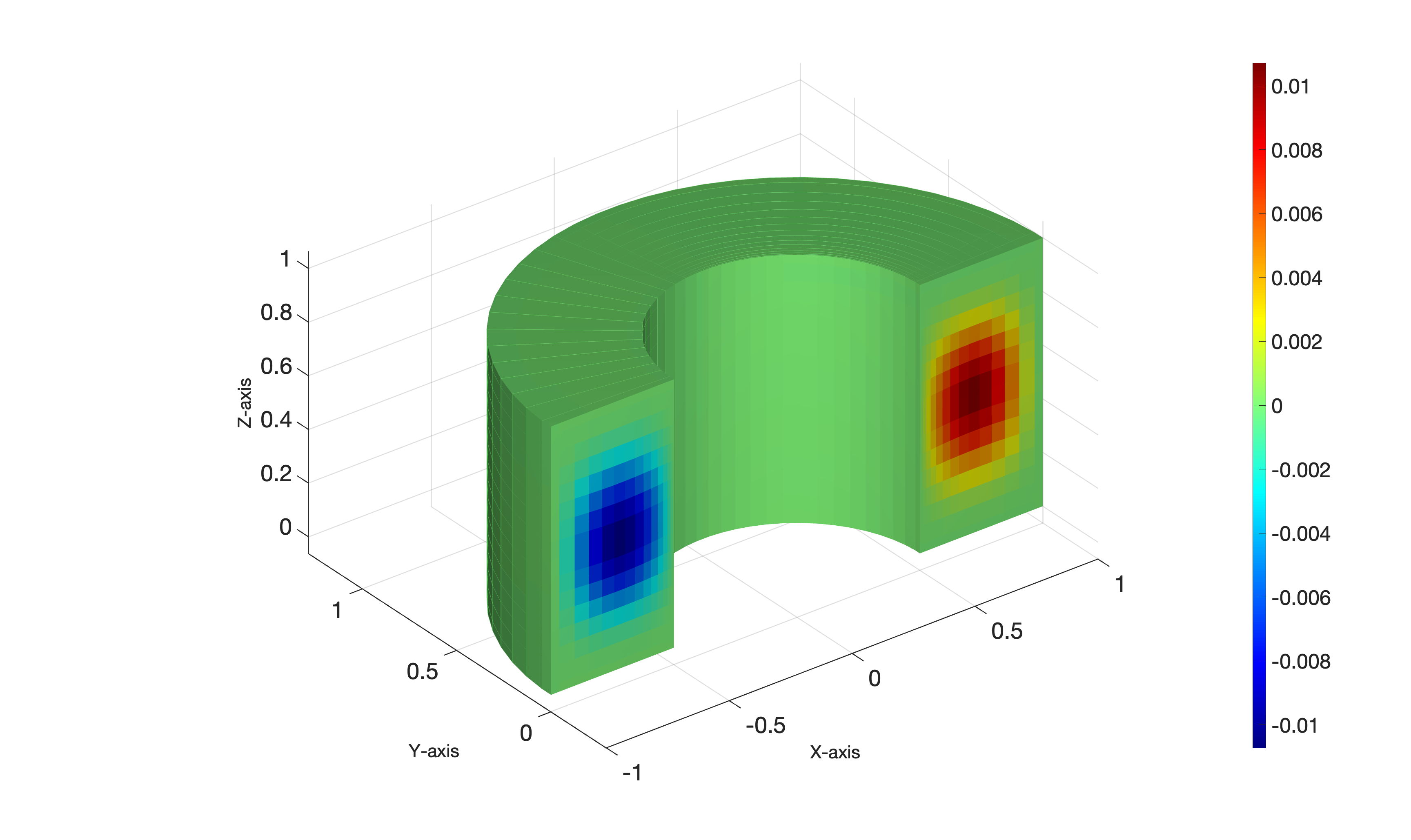}
\caption*{(d) $z$-displacement}
\endminipage
\caption{Mesh grid and displacement fields of the cylindrical structure. Because a zero boundary condition is applied to the outer surface of the ring, half of the structure is omitted from this visualization.}
\label{Fig.Ring-displacement}
\end{center}
\end{figure}

\begin{figure}
\begin{center}
\minipage{0.5\textwidth}
\includegraphics[width=0.99\textwidth]{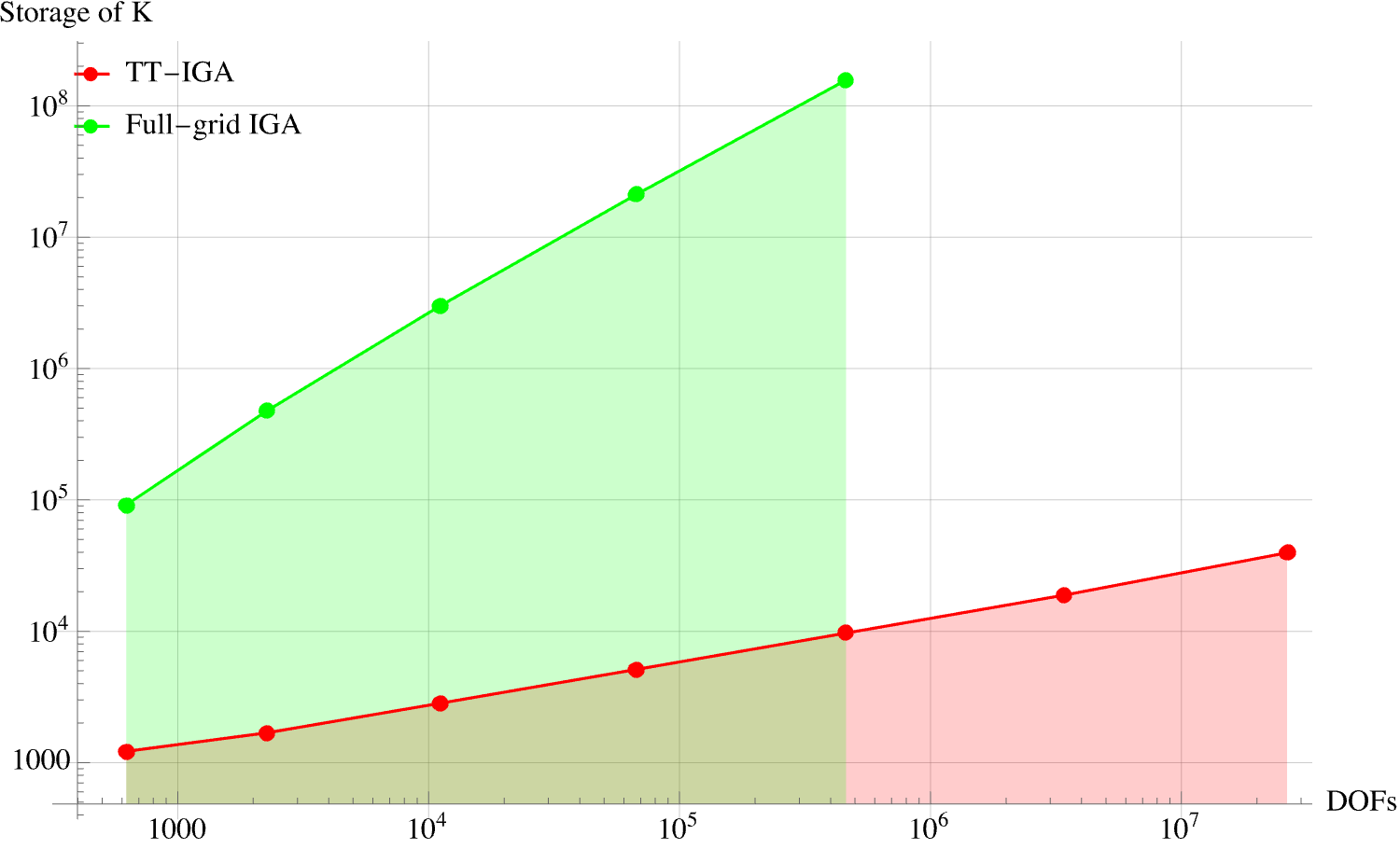}
\caption*{(a) Stiffness matrix}
\endminipage
\hfill
\minipage{0.5\textwidth}
\includegraphics[width=0.99\textwidth]{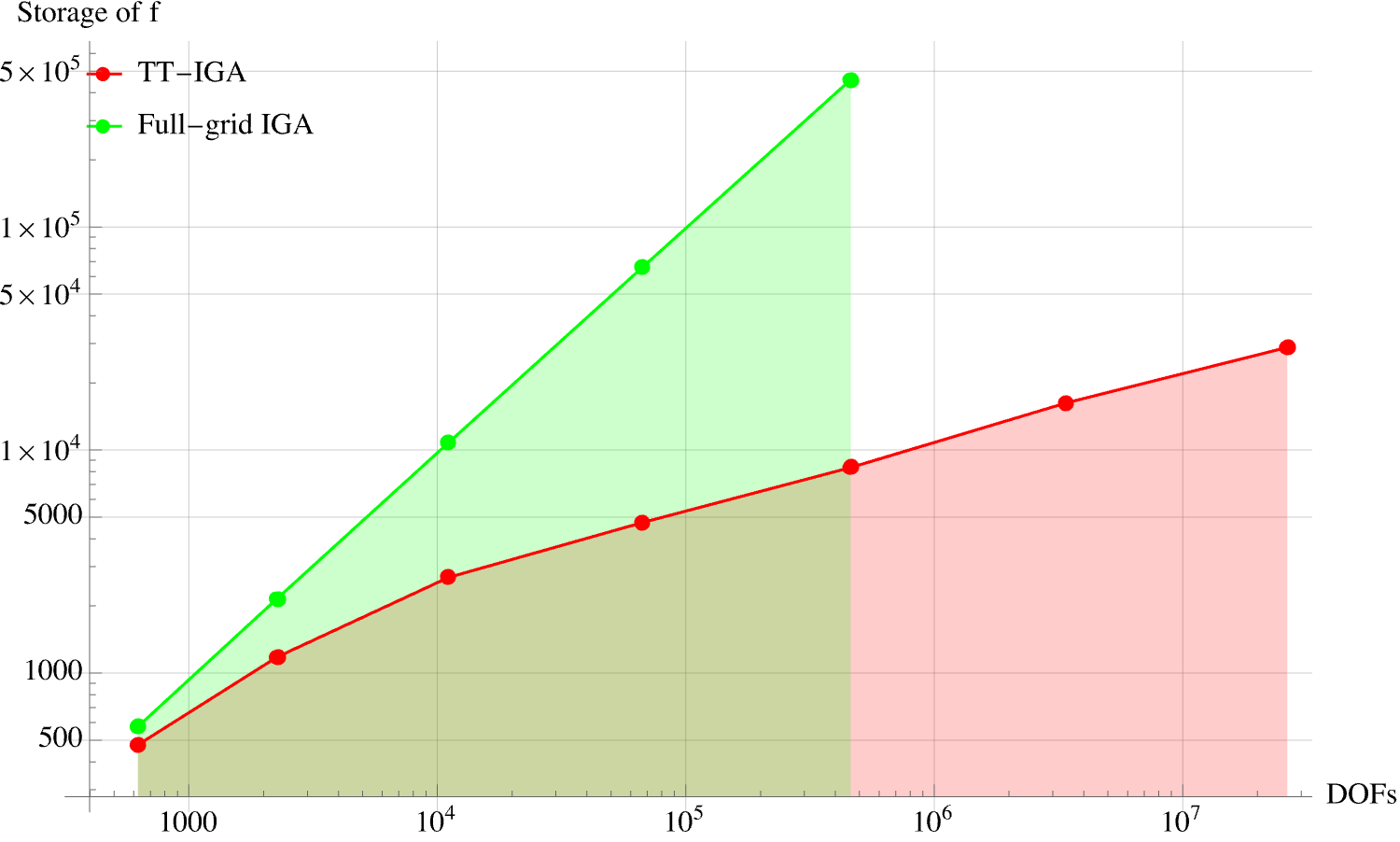}
\caption*{(b) Storage of applied force $\boldsymbol{f}$}
\endminipage
\hfill
\minipage{0.5\textwidth}
\includegraphics[width=0.99\textwidth]{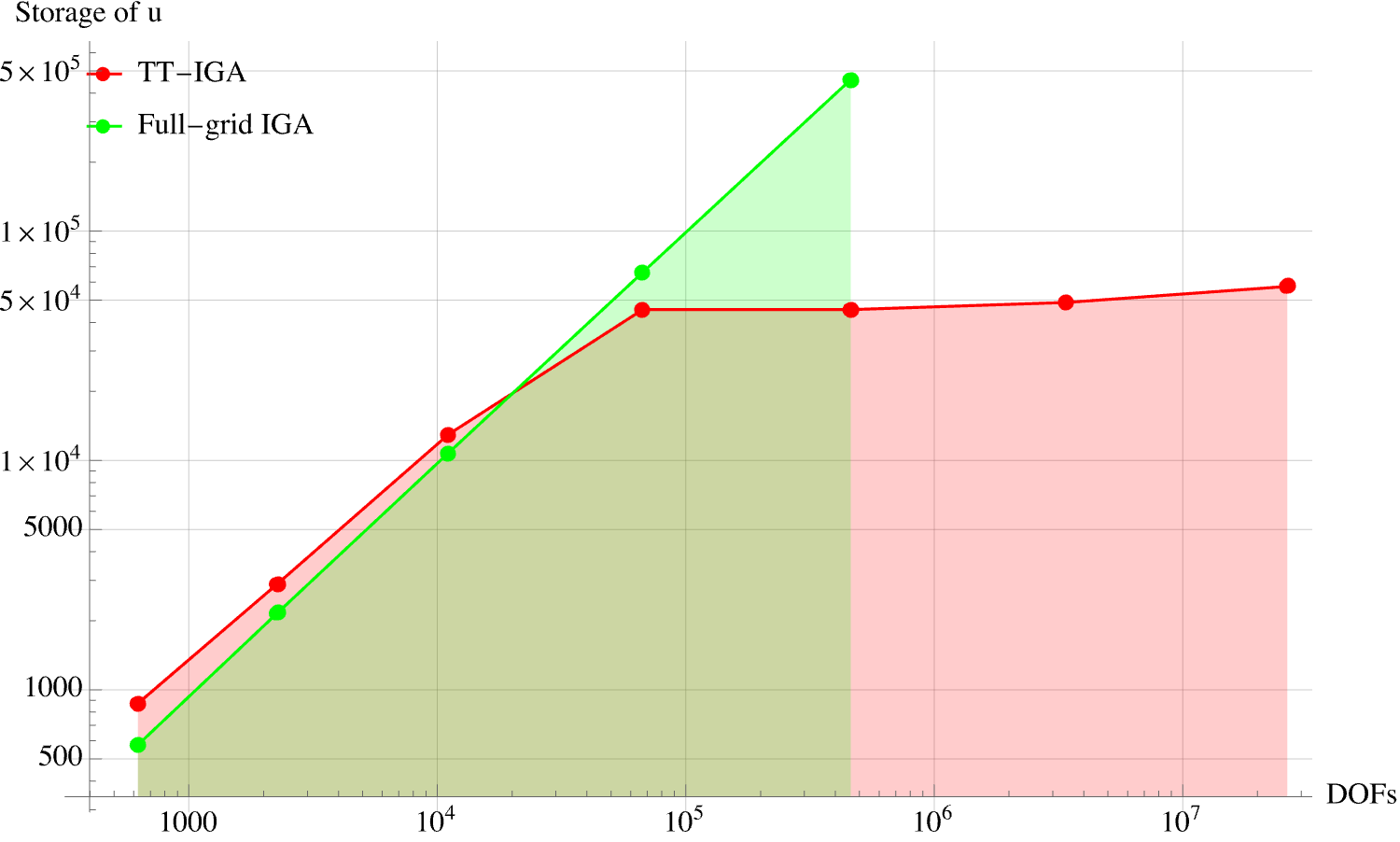}
\caption*{(c) Displacement}
\endminipage
\hfill
\minipage{0.5\textwidth}
\includegraphics[width=0.99\textwidth]{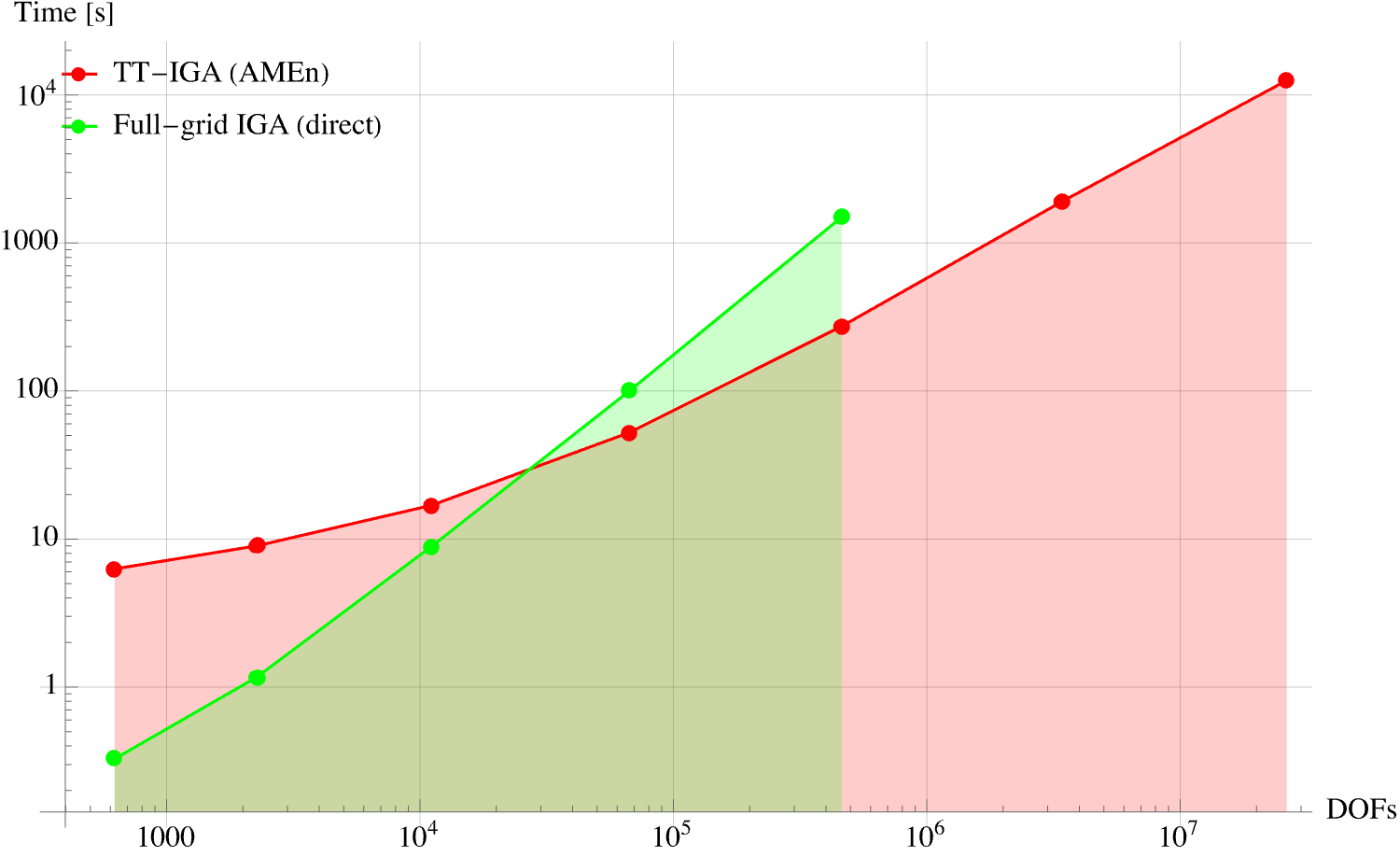}
\caption*{(d) Total computational time}
\endminipage
\hfill
\minipage{0.5\textwidth}
\includegraphics[width=0.99\textwidth]{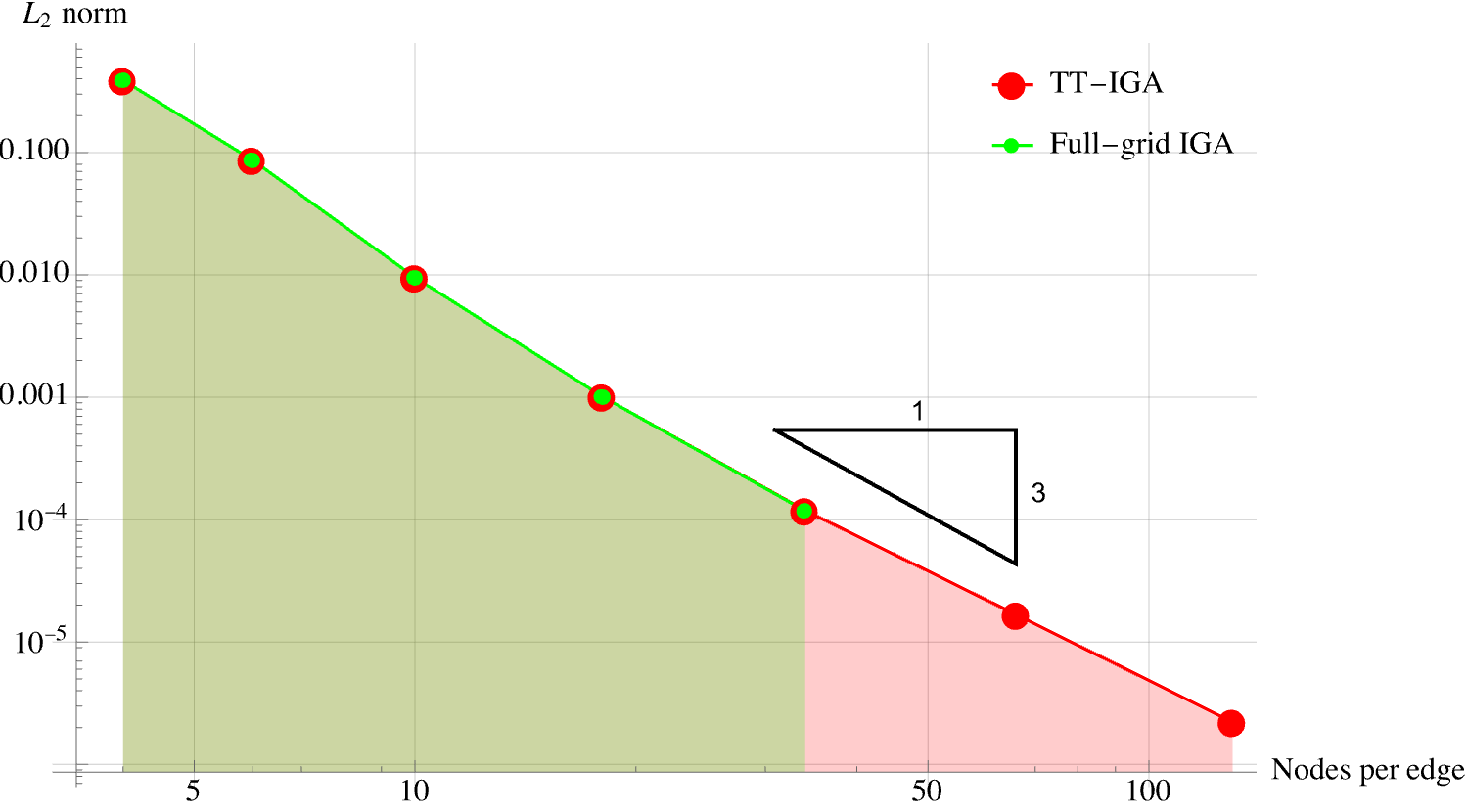}
\caption*{(e) $L_2$ norm}
\endminipage
\caption{Numerical results for the 3D cylindrical structure. Storage costs are reported for (a) the stiffness matrix \(\boldsymbol{K}\), (b) the applied force vector \(\boldsymbol{f}\), and (c) the displacement solution \(\boldsymbol{u}\). The solution time and \(L_2\)-error convergence are shown in (d) and (e), respectively.}
\label{Fig.Ring-results}
\end{center}
\end{figure}

\Figref{Fig.Ring-displacement} shows the mesh with a \(69 \times 18 \times 18\) discretization and the displacement fields obtained from the numerical solution. In this verification example, the TT formulation solved with AMEn is compared directly with the conventional full-grid IGA formulation. \Figref{Fig.Ring-results} reports the storage of the stiffness matrix \(\boldsymbol{K}\), force vector \(\boldsymbol{f}\), and displacement solution \(\boldsymbol{u}\), together with the elapsed solution time and the \(L_2\)-error convergence.

The operator and force-vector results show that the TT representation requires substantially fewer stored entries than the assembled sparse full-grid representation, and the gap increases under refinement. The displacement vector exhibits a crossover behavior: for the smallest systems, the TT representation carries additional rank-related overhead, whereas at larger problem sizes the displacement becomes increasingly more compact than its full-grid counterpart. A similar crossover is observed in the solution time. The sparse direct solver is faster for the smallest systems, while the AMEn-based TT solution becomes more favorable as the problem size increases and continues to larger discretizations after the full-grid calculations become impractical.

Most importantly, \Figref{Fig.Ring-results}e shows that the TT and full-grid solutions exhibit essentially the same \(L_2\)-error convergence, with approximately third-order behavior for the quadratic basis used in this example. Thus, the compression introduced by the TT representation does not alter the observed discretization accuracy. This example therefore establishes the numerical consistency of the proposed TT-IGA elasticity formulation before its scalability is examined on the larger geometries below.

\subsection{Open hemisphere}
\begin{figure}
\begin{center}
\includegraphics[width=0.60\textwidth]{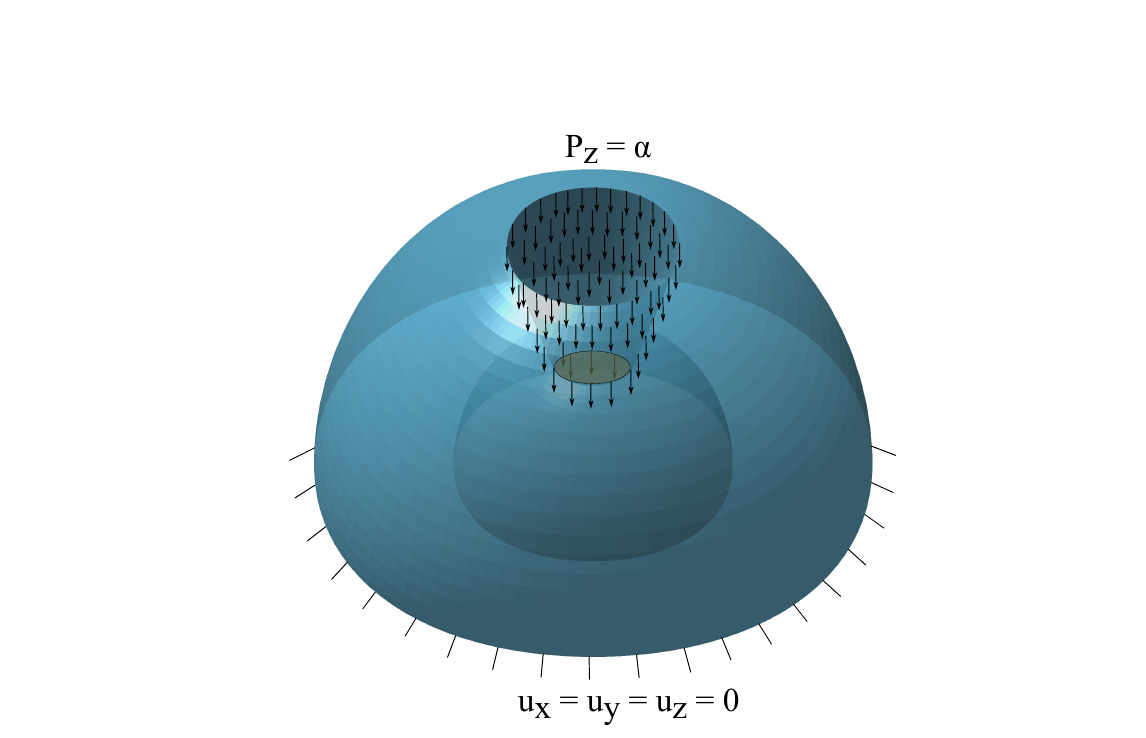}
\end{center}
\caption{Open hemispherical structure clamped at the bottom surface and subjected to a distributed surface force on the inner surface.}
\label{Fig.Hemi}
\end{figure}

\begin{figure}
\begin{center}
\minipage{0.5\textwidth}
\includegraphics[width=0.70\textwidth]{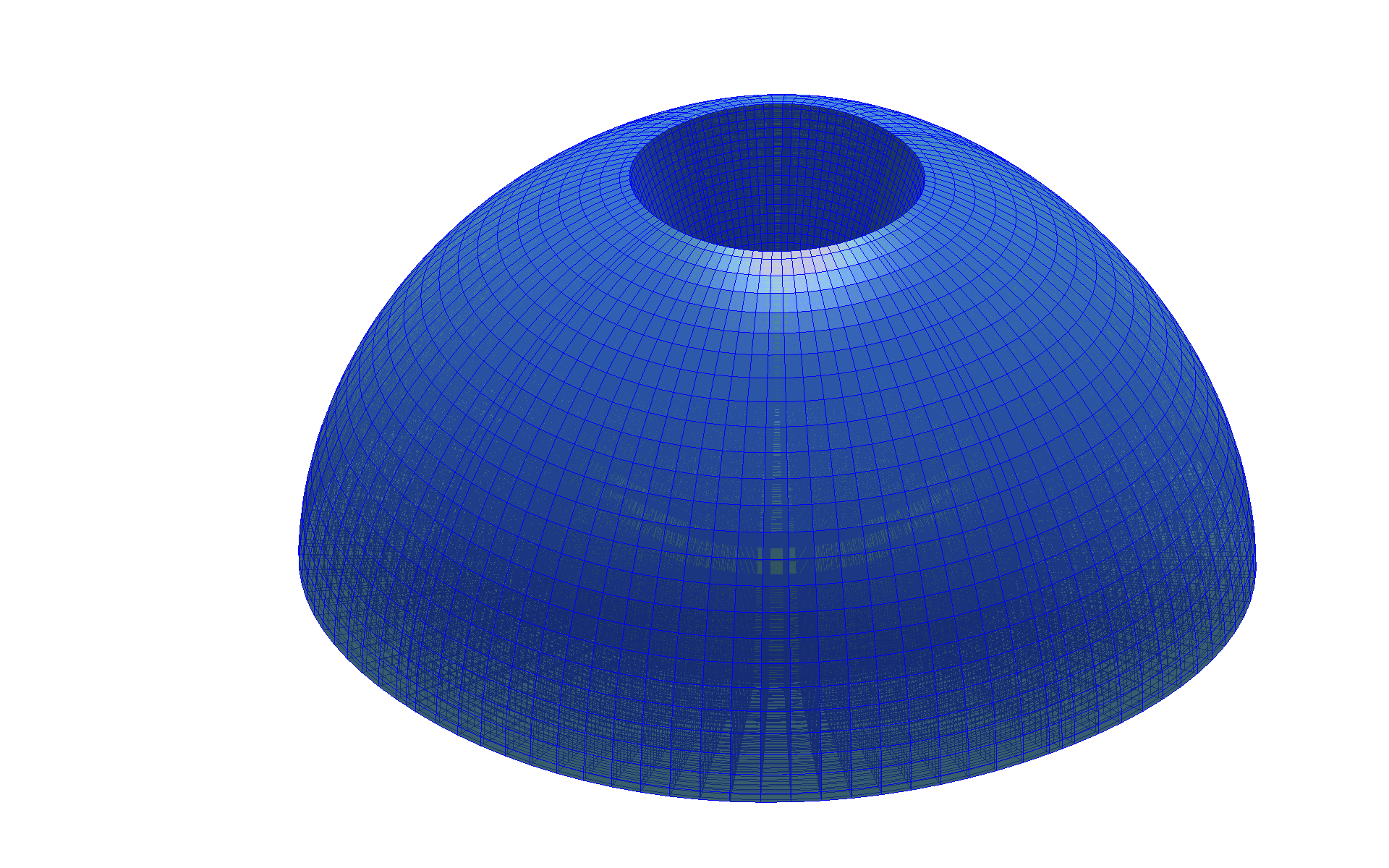}
\caption*{(a) Control mesh grid $101\times 26 \times 26$}
\endminipage
\hfill
\minipage{0.5\textwidth}
\includegraphics[width=1.00\textwidth]{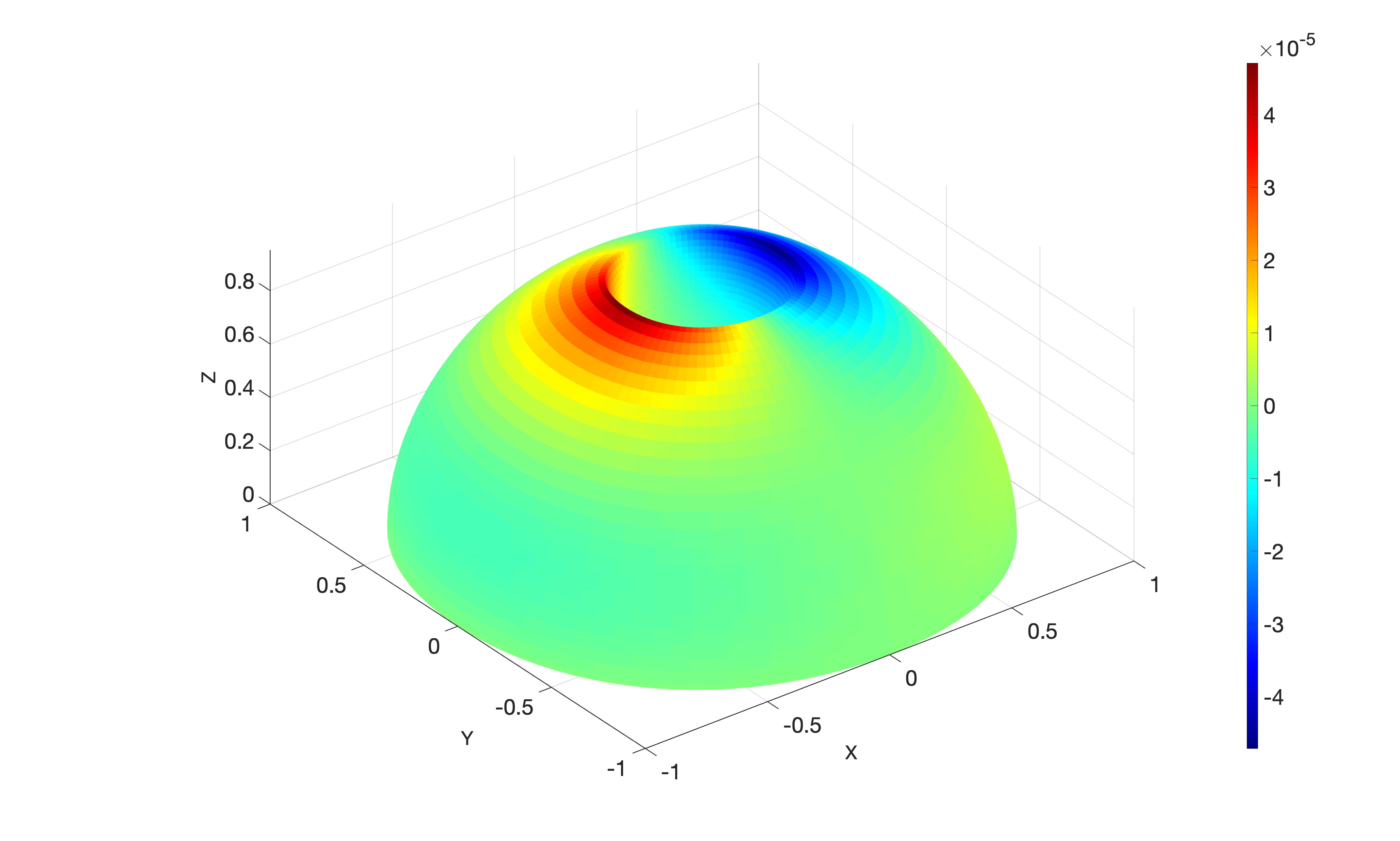}
\caption*{(b) $x$-displacement}
\endminipage
\hfill
\minipage{0.5\textwidth}
\includegraphics[width=1.00\textwidth]{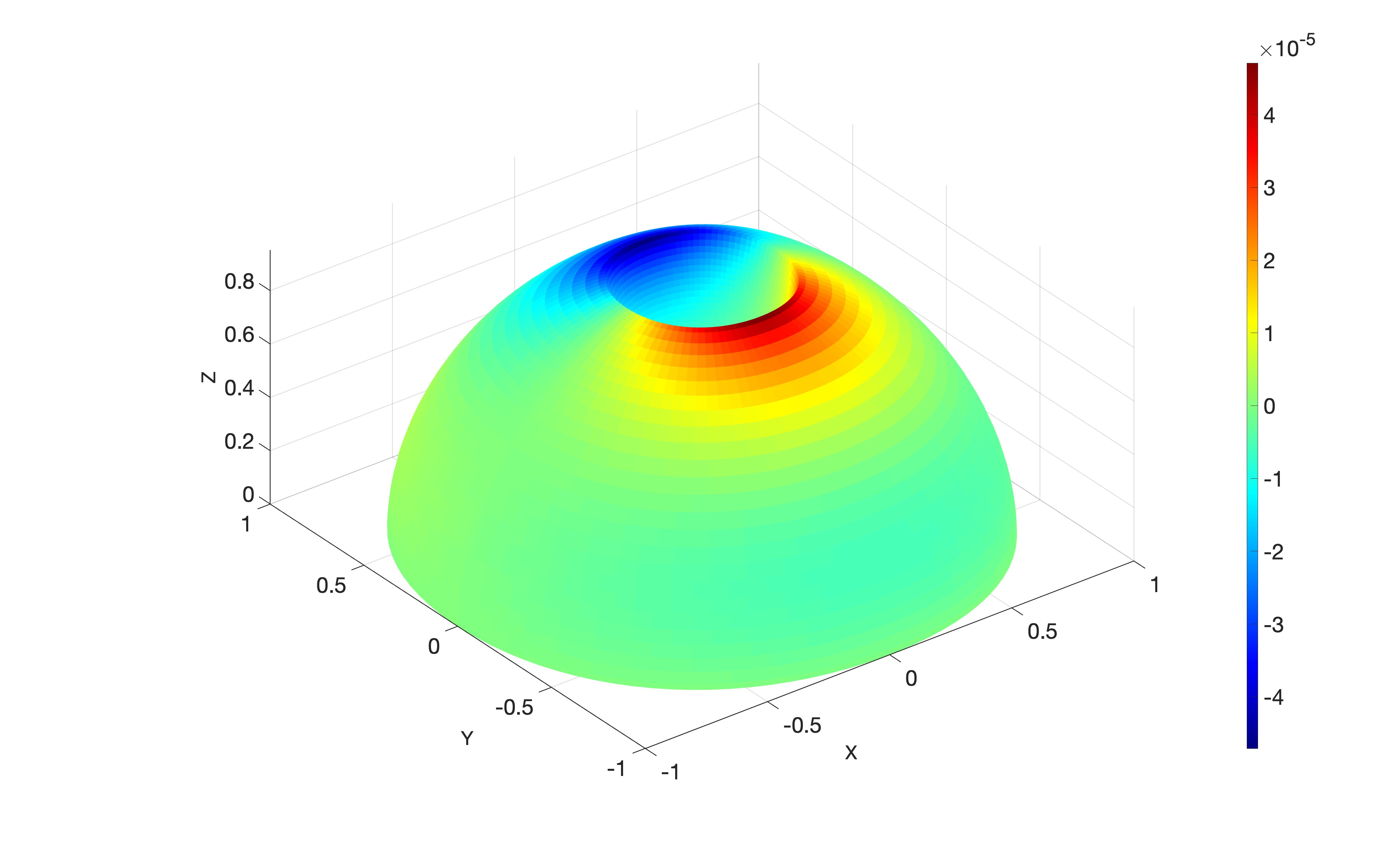}
\caption*{(c) $y$-displacement}
\endminipage
\hfill
\minipage{0.5\textwidth}
\includegraphics[width=1.00\textwidth]{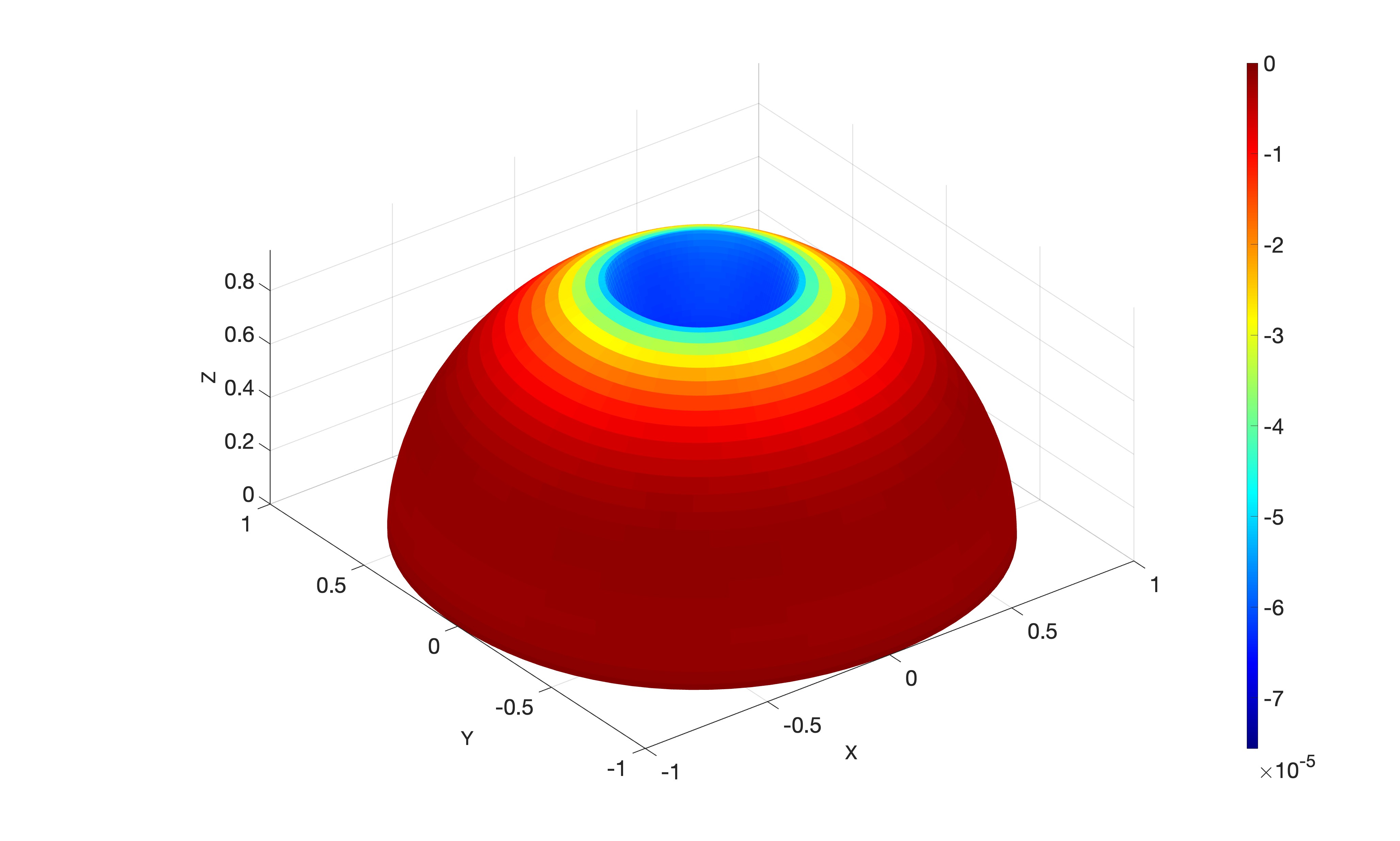}
\caption*{(d) $z$-displacement}
\endminipage
\caption{Mesh and displacement fields of the open hemispherical structure.}
\label{Fig.Hemi-displacement}
\end{center}
\end{figure}

In this numerical example, an open hemispherical shell subjected to a clamped boundary condition at its bottom edge and a distributed load on the inner surface is examined, as shown in \Figref{Fig.Hemi}. The geometry is defined by the inner radius \(R_{\mathrm{in}}=0.5\,\mathrm{m}\), the outer radius \(R_{\mathrm{out}}=1.0\,\mathrm{m}\), and an \(18^\circ\) opening at the top. Quadratic basis functions, \(p_1=p_2=p_3=2\), are used. A Neumann boundary condition is applied on the inner surface with \(P_z=-100\,\mathrm{N/m^2}\). The displacement field obtained using a \(101\times 26\times 26\) control mesh is presented in \Figref{Fig.Hemi-displacement}.

\begin{figure}
\begin{center}
\minipage{0.5\textwidth}
\includegraphics[width=0.99\textwidth]{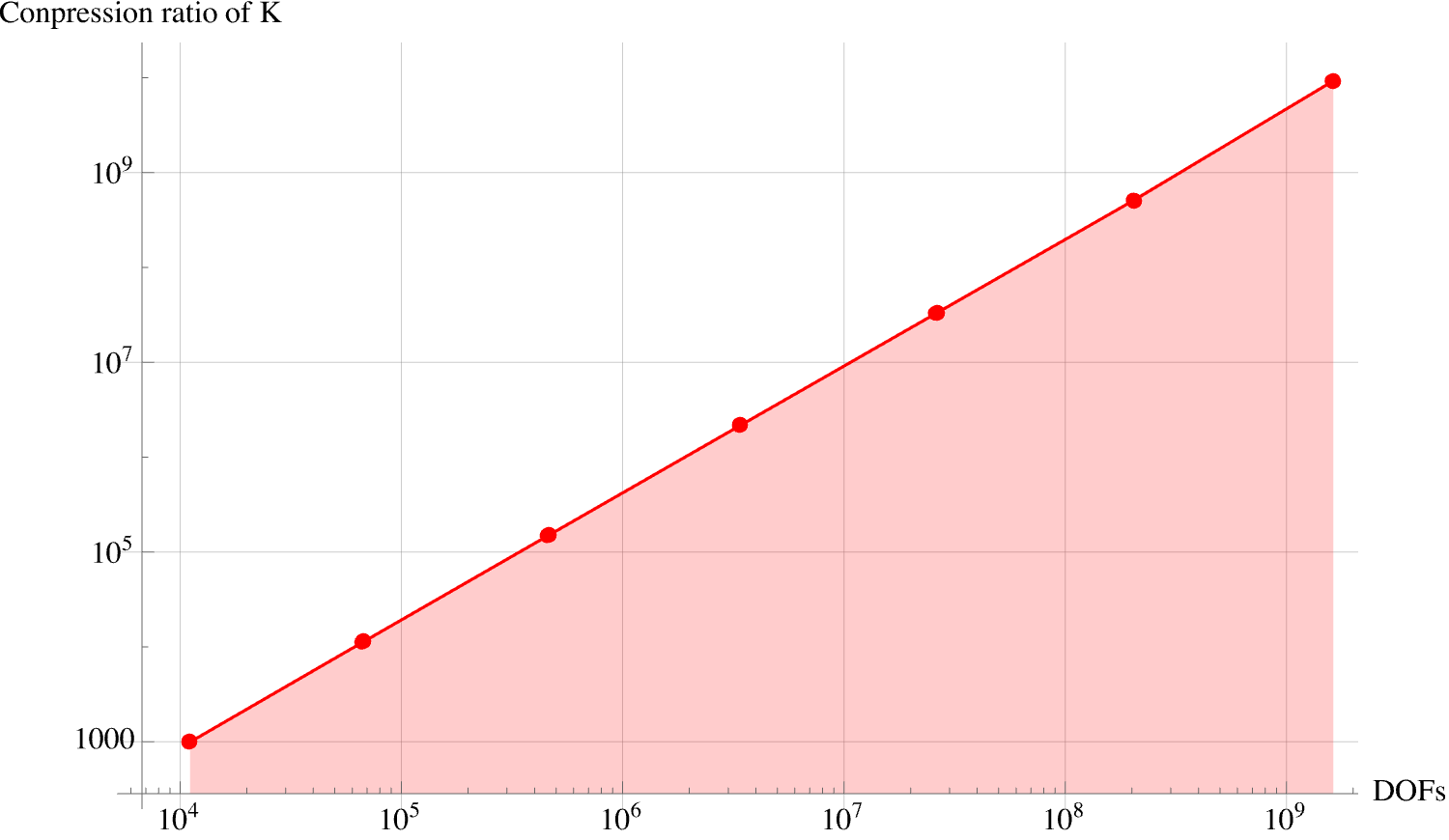}
\caption*{(a) Stiffness matrix}
\endminipage
\hfill
\minipage{0.5\textwidth}
\includegraphics[width=0.99\textwidth]{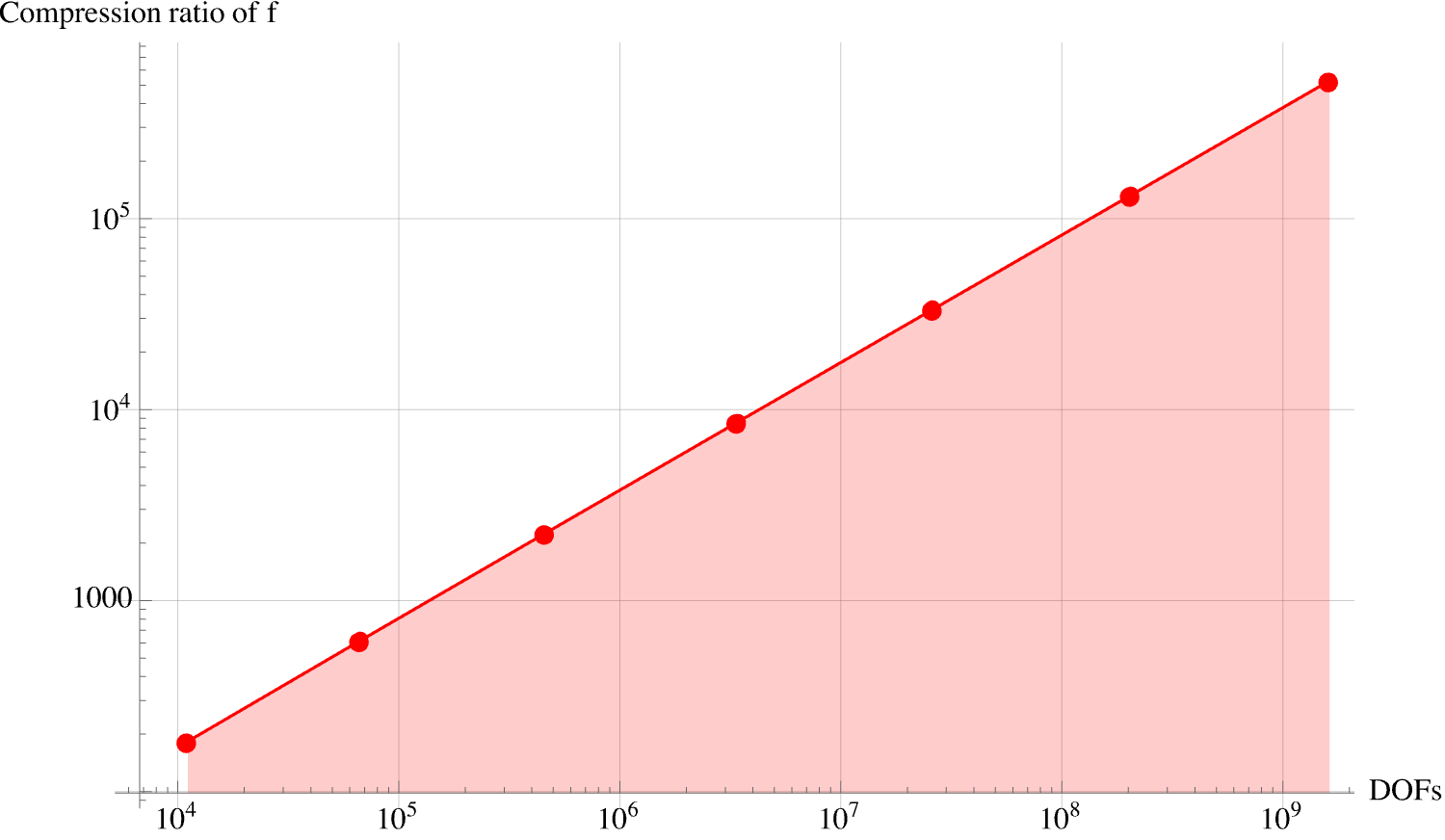}
\caption*{(b) Applied force $\boldsymbol{f}$}
\endminipage
\hfill
\minipage{0.5\textwidth}
\includegraphics[width=0.99\textwidth]{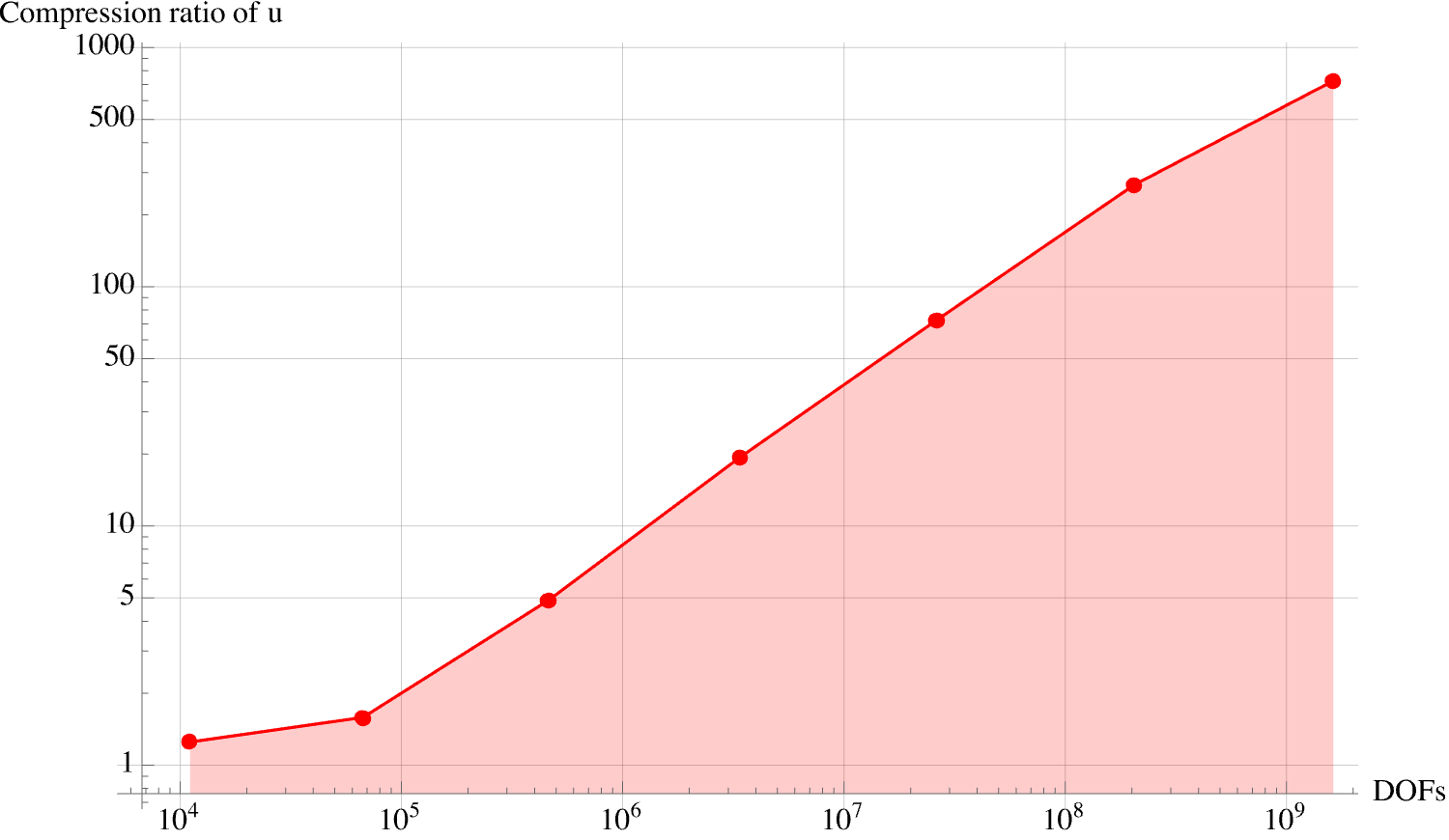}
\caption*{(c) Displacement}
\endminipage
\hfill
\minipage{0.5\textwidth}
\includegraphics[width=0.99\textwidth]{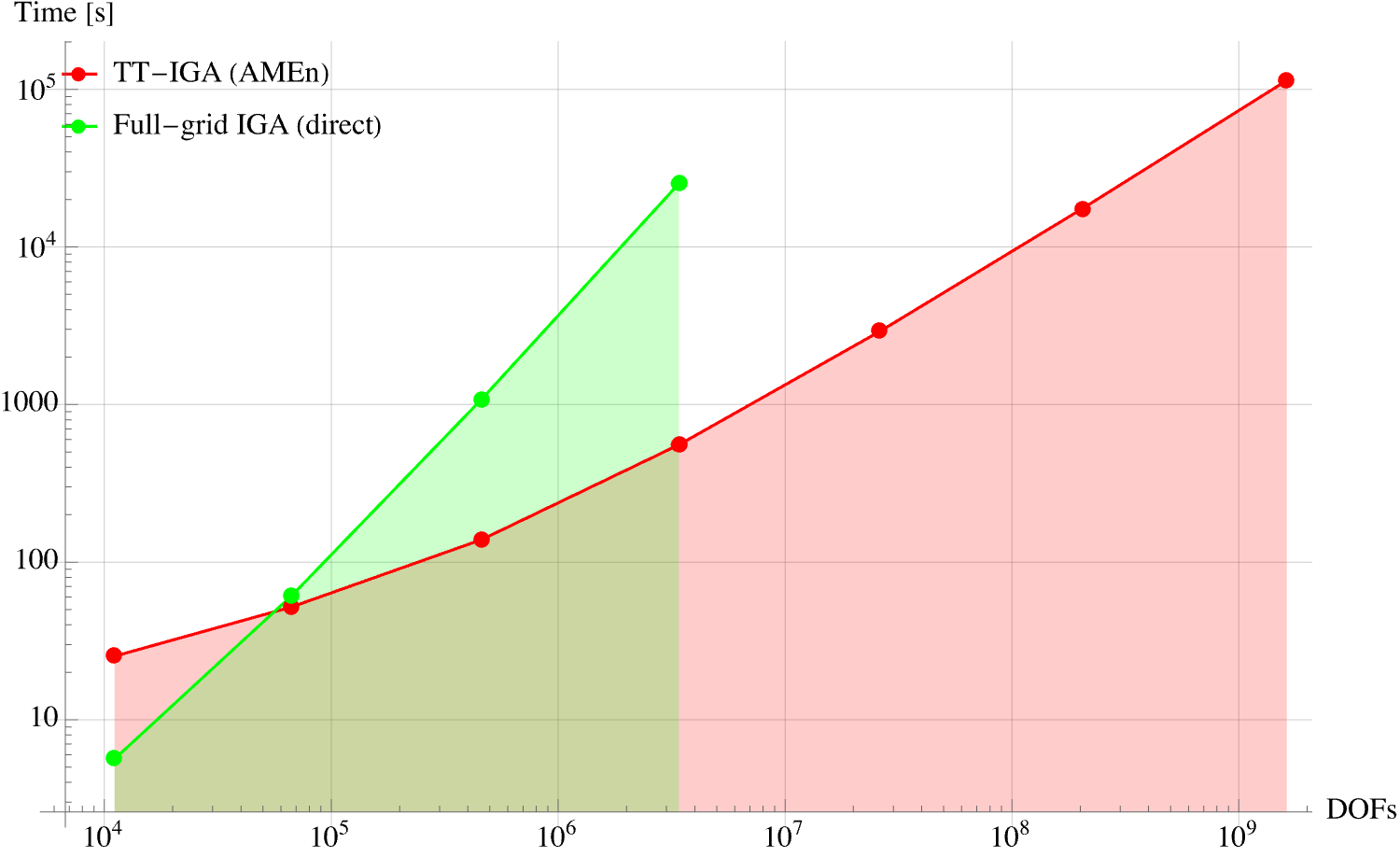}
\caption*{(d) Total computational time}
\endminipage
\caption{Compression ratios of the stiffness matrix $\boldsymbol{K}$, force vector $\boldsymbol{f}$, and displacement solution $\boldsymbol{u}$, together with TT-AMEn and full-grid computational times for the open hemispherical geometry.}
\label{Fig.Hemi-compression}
\end{center}
\end{figure}

\Figref{Fig.Hemi-compression} shows the compression ratios of the stiffness operator \(\boldsymbol{K}\), force vector \(\boldsymbol{f}\), and displacement solution \(\boldsymbol{u}\), together with the TT-AMEn and full-grid computational times. The compression ratios increase strongly with refinement, indicating that the storage required by the TT cores grows much more slowly than the corresponding full-grid representation. The stiffness operator exhibits particularly strong compression, while the displacement compression also increases steadily as the number of degrees of freedom grows. The constant surface traction remains highly compressible because its construction involves low-rank spline quantities and the surface Jacobian rather than a general three-dimensional sampled field.

The timing comparison shows the expected overhead of the TT representation at the smallest discretization, but this behavior reverses as the problem size increases. Over the range where both approaches can be executed, the total computational cost of the full-grid workflow grows much more rapidly, whereas the TT-AMEn workflow retains a substantially milder growth rate. The TT formulation then continues to problem sizes for which explicit full-grid assembly and solution are no longer practical. The largest open-hemisphere calculation reaches approximately \(1.6\times10^9\) degrees of freedom, illustrating that the principal advantage of the formulation is not merely compression of an already tractable system, but the ability to represent and solve structured three-dimensional elasticity problems beyond the accessible full-grid range.
\subsection{Square plate with a circular hole}
\begin{figure}
\begin{center}
\includegraphics[width=0.70\textwidth]{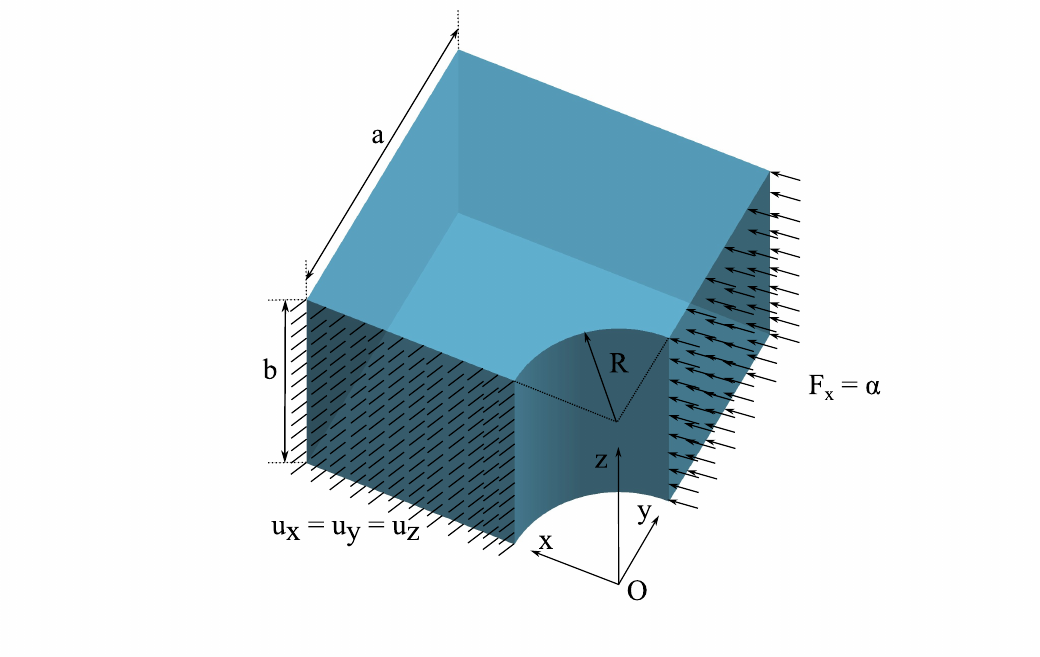}
\end{center}
\caption{Square plate with a circular hole: the structure is clamped at the face $y = 0$ and a constant traction force $F_x = \alpha$ is applied on the face $x = 0$.}
\label{Fig.SquareHole}
\end{figure}

\begin{figure}
\begin{center}
\minipage{0.5\textwidth}
\includegraphics[width=0.65\textwidth]{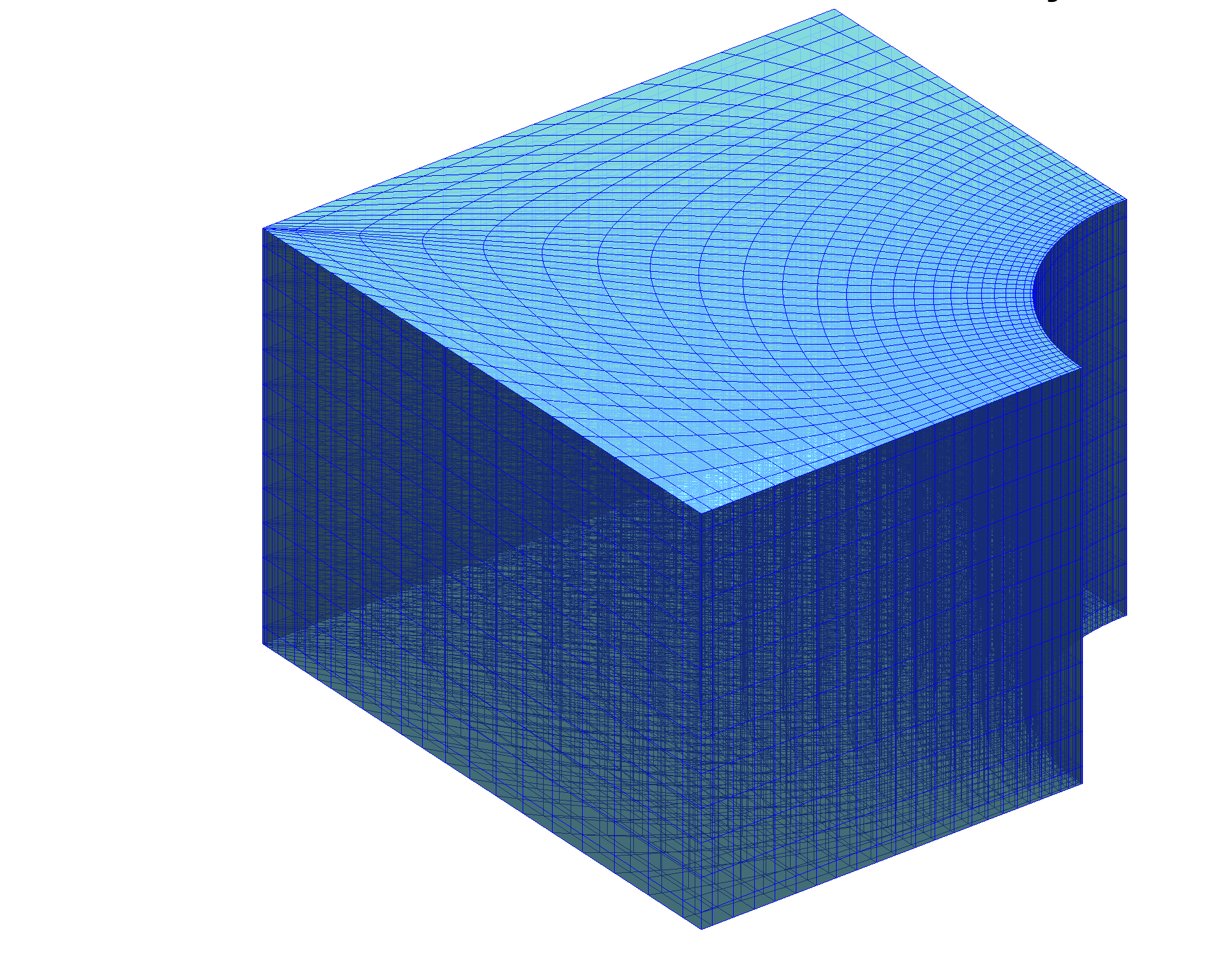}
\caption*{(a) Control mesh grid $50 \times 26 \times 14$}
\endminipage
\minipage{0.5\textwidth}
\includegraphics[width=1.0\textwidth]{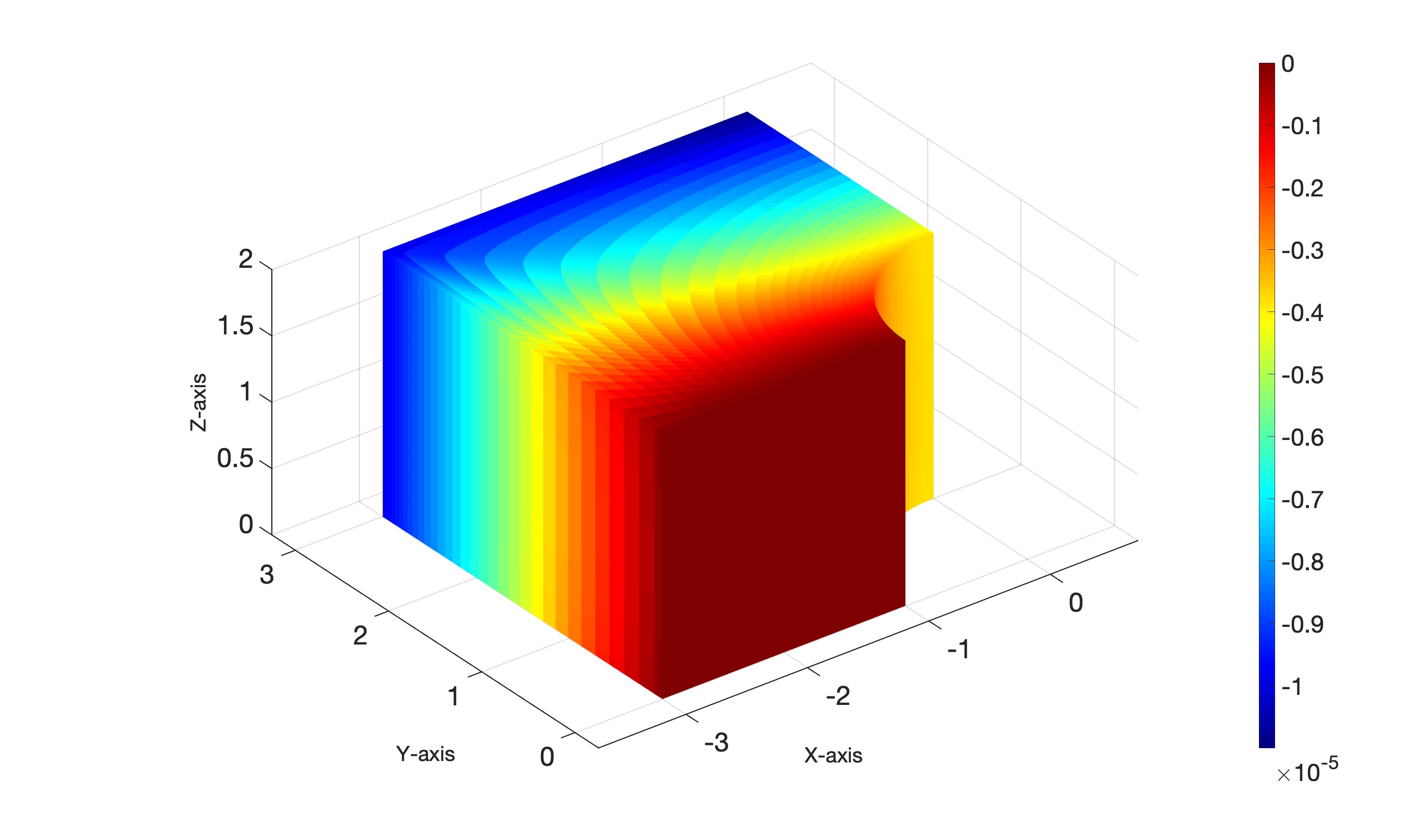}
\caption*{(b) $x$-displacement}
\endminipage
\hfill
\minipage{0.5\textwidth}
\includegraphics[width=1.0\textwidth]{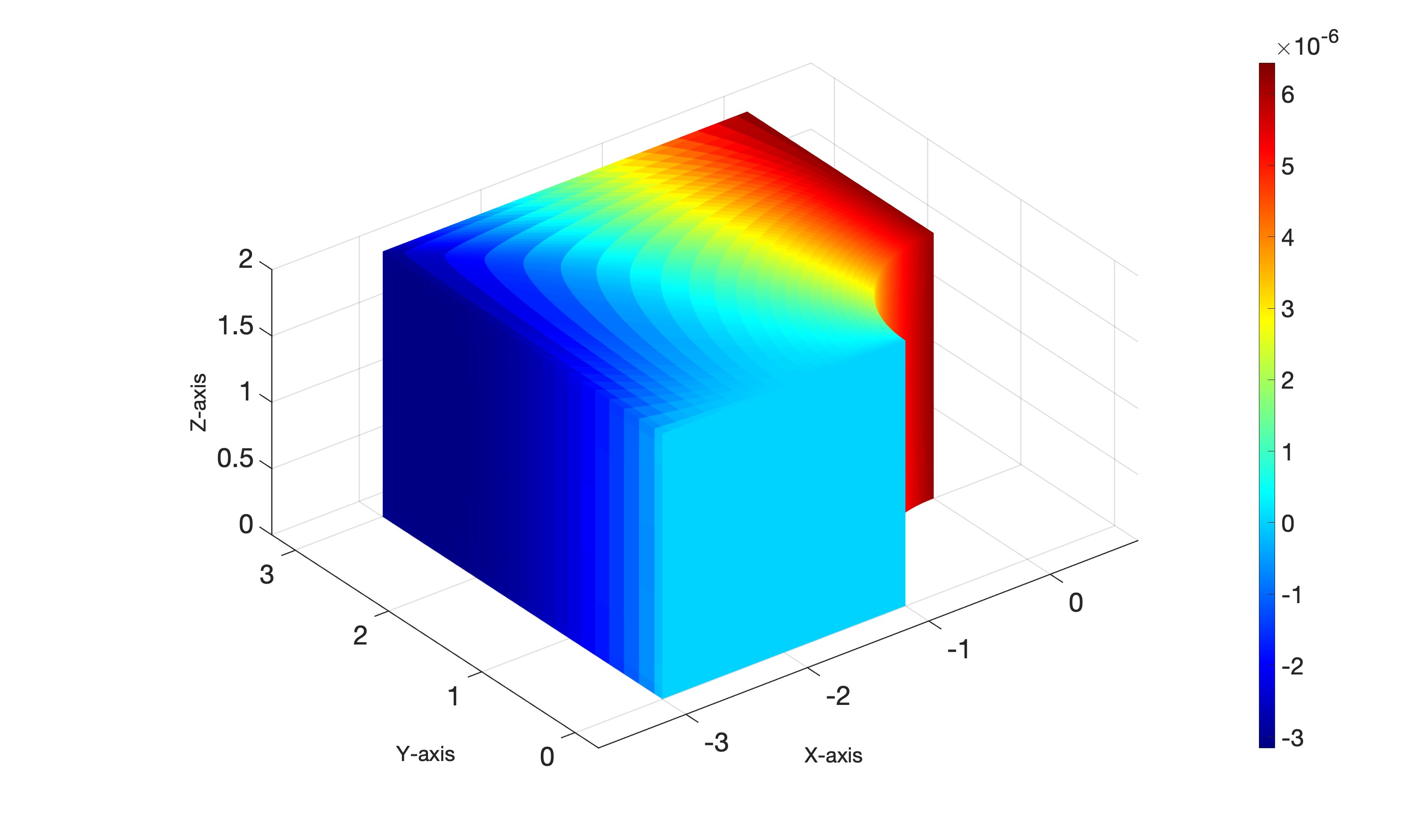}
\caption*{(c) $y$-displacement}
\endminipage
\hfill
\minipage{0.5\textwidth}
\includegraphics[width=1.0\textwidth]{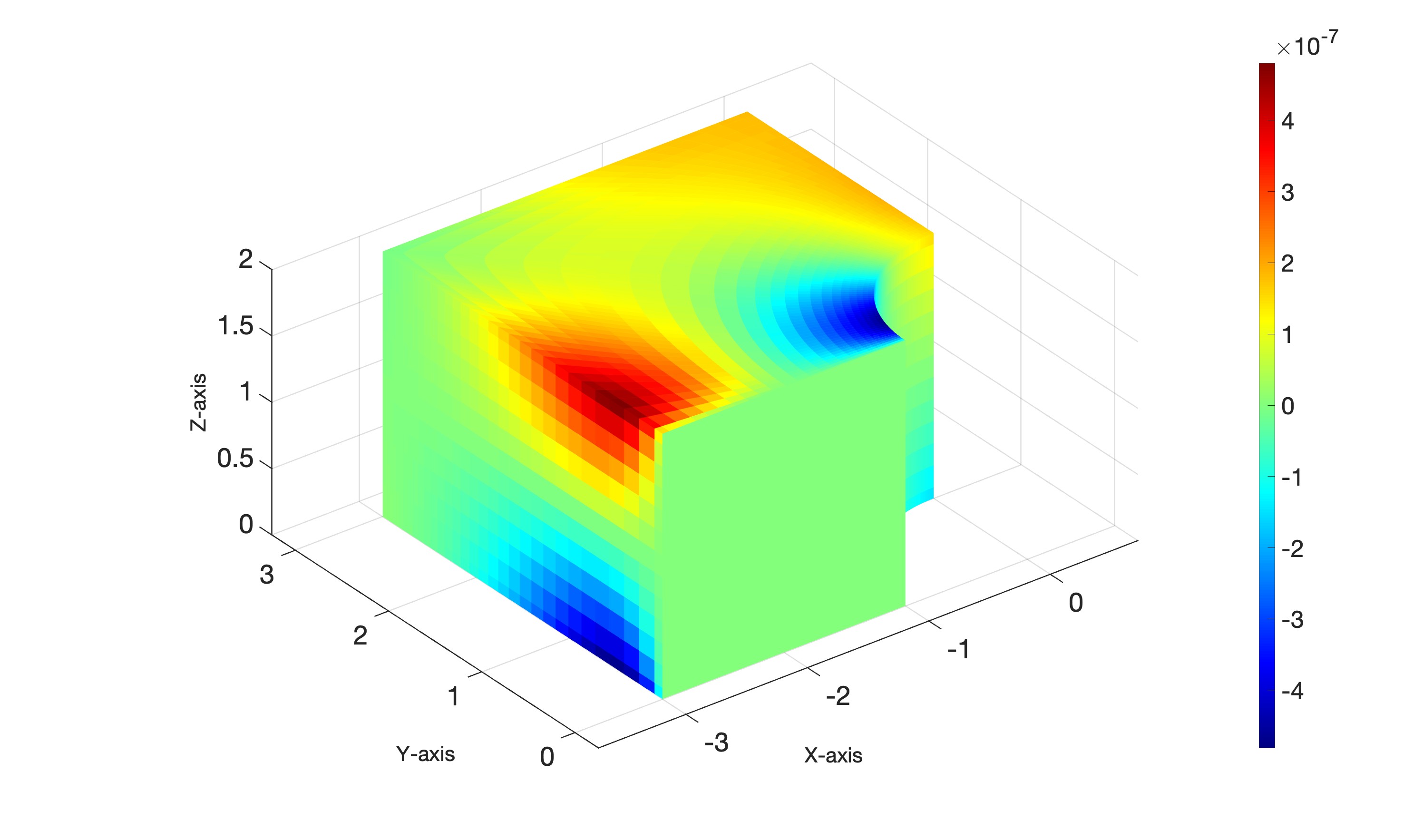}
\caption*{(d) $z$-displacement}
\endminipage
\caption{Mesh and displacement fields of the square plate with a circular hole.}
\label{Fig.square-hole-displacement}
\end{center}
\end{figure}

A three-dimensional square plate with a circular hole, as illustrated in \Figref{Fig.SquareHole}, is investigated in this example. The geometric dimensions are \(R=1\,\mathrm{m}\), \(b=2\,\mathrm{m}\), and \(a=3\,\mathrm{m}\). Quadratic basis functions are used in all three parametric directions, i.e., \(p_1=p_2=p_3=2\). Homogeneous Dirichlet boundary conditions are imposed on the plane \(y=0\), and a constant surface traction \(F_x=\alpha=-10\,\mathrm{N/m^2}\) is applied on the plane \(x=0\). The resulting displacement fields are shown in \Figref{Fig.square-hole-displacement} for a control mesh of \(50 \times 26 \times 14\).

\begin{figure}
\begin{center}
\minipage{0.5\textwidth}
\includegraphics[width=0.99\textwidth]{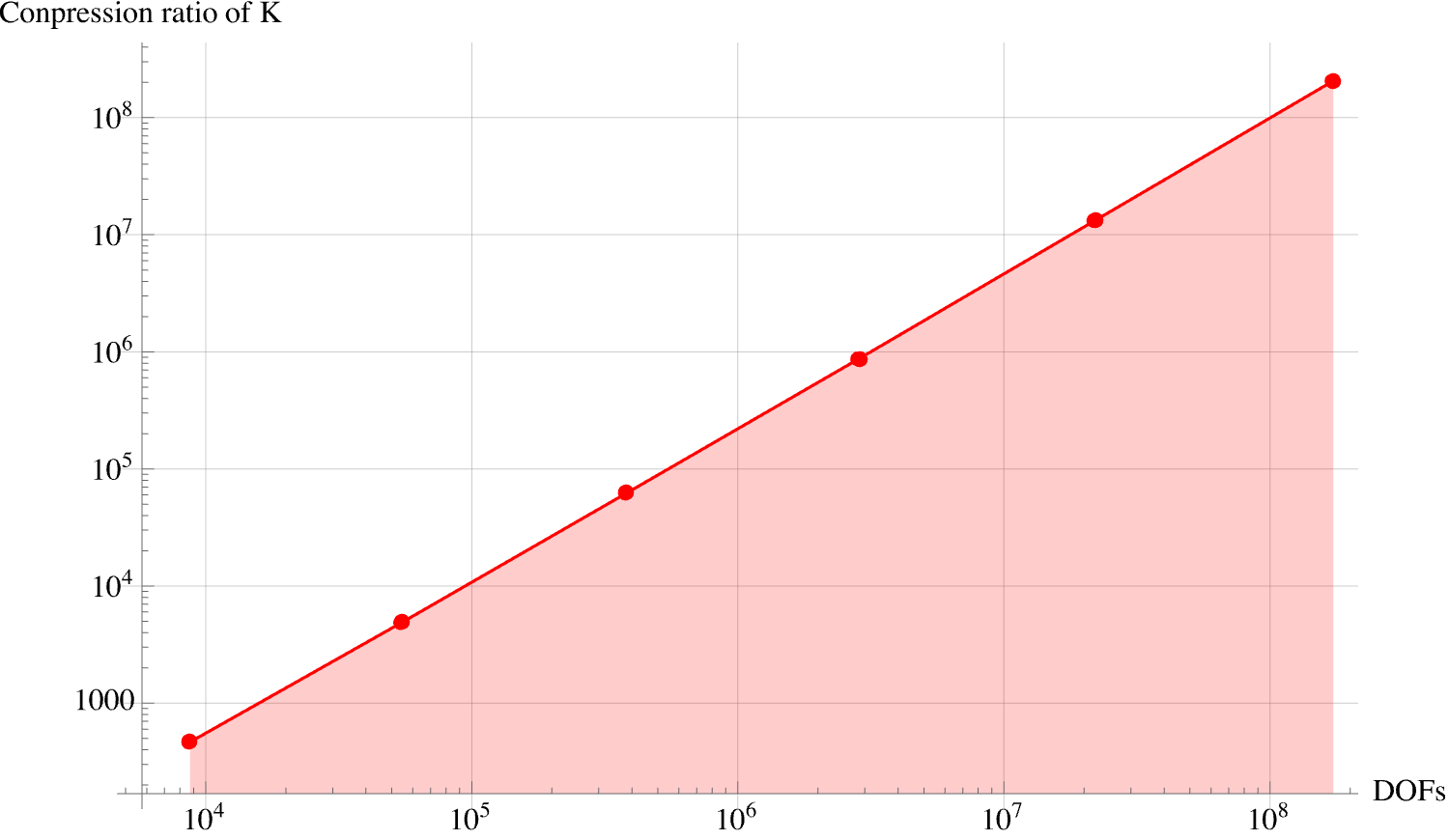}
\caption*{(a) Stiffness matrix}
\endminipage
\hfill
\minipage{0.5\textwidth}
\includegraphics[width=0.99\textwidth]{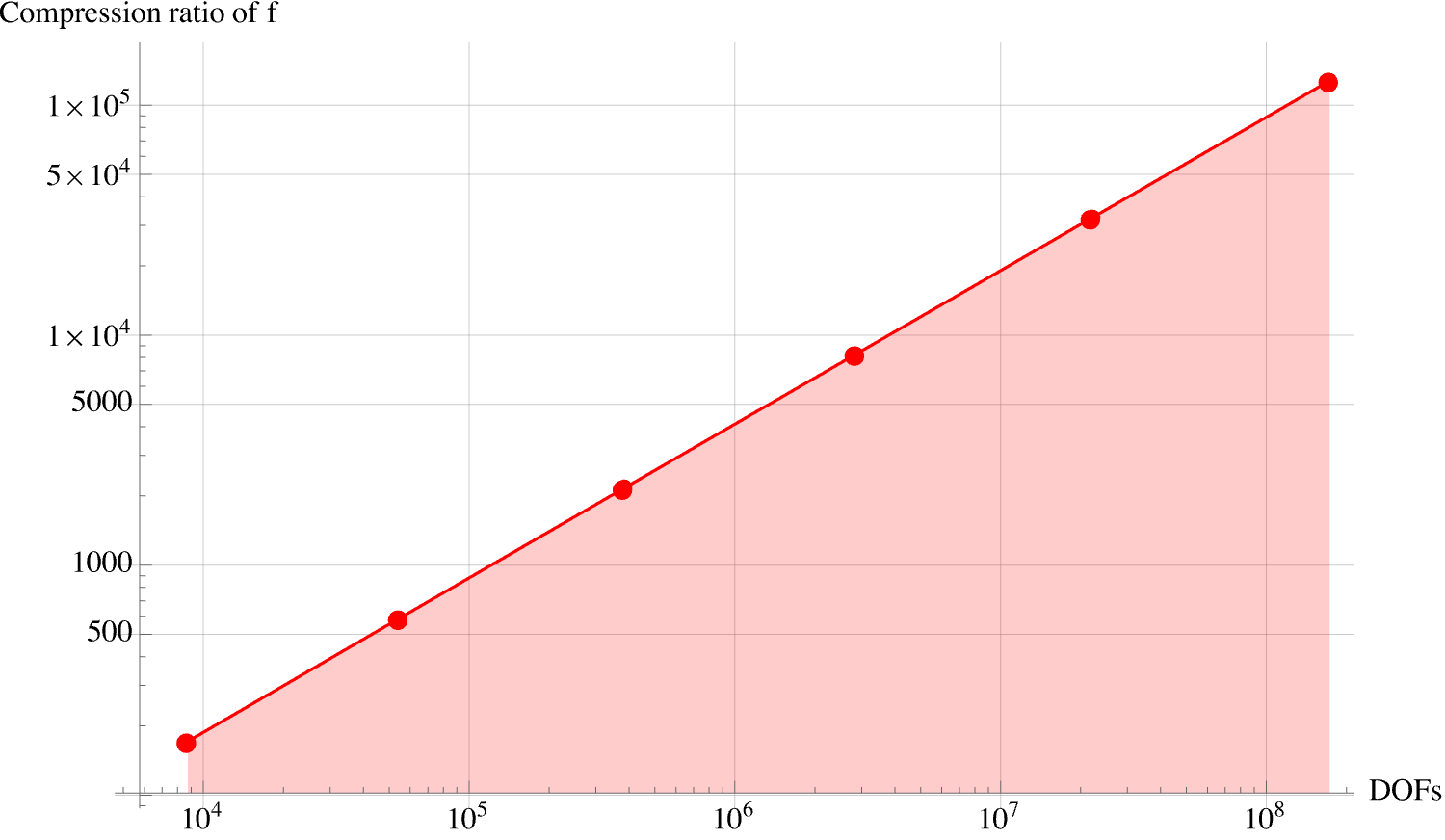}
\caption*{(b) Applied force $\boldsymbol{f}$}
\endminipage
\hfill
\minipage{0.5\textwidth}
\includegraphics[width=0.99\textwidth]{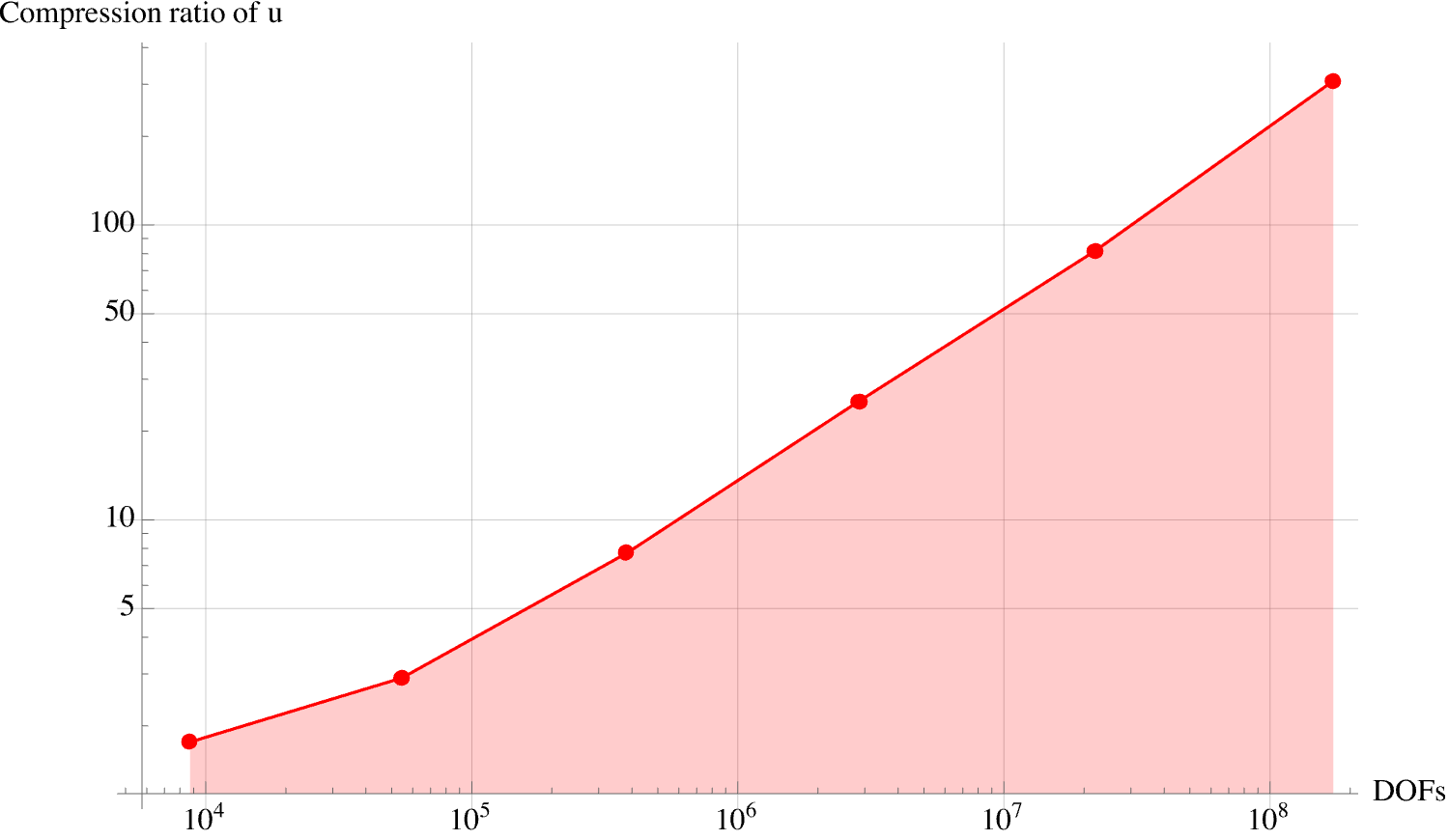}
\caption*{(c) Displacement}
\endminipage
\hfill
\minipage{0.5\textwidth}
\includegraphics[width=0.99\textwidth]{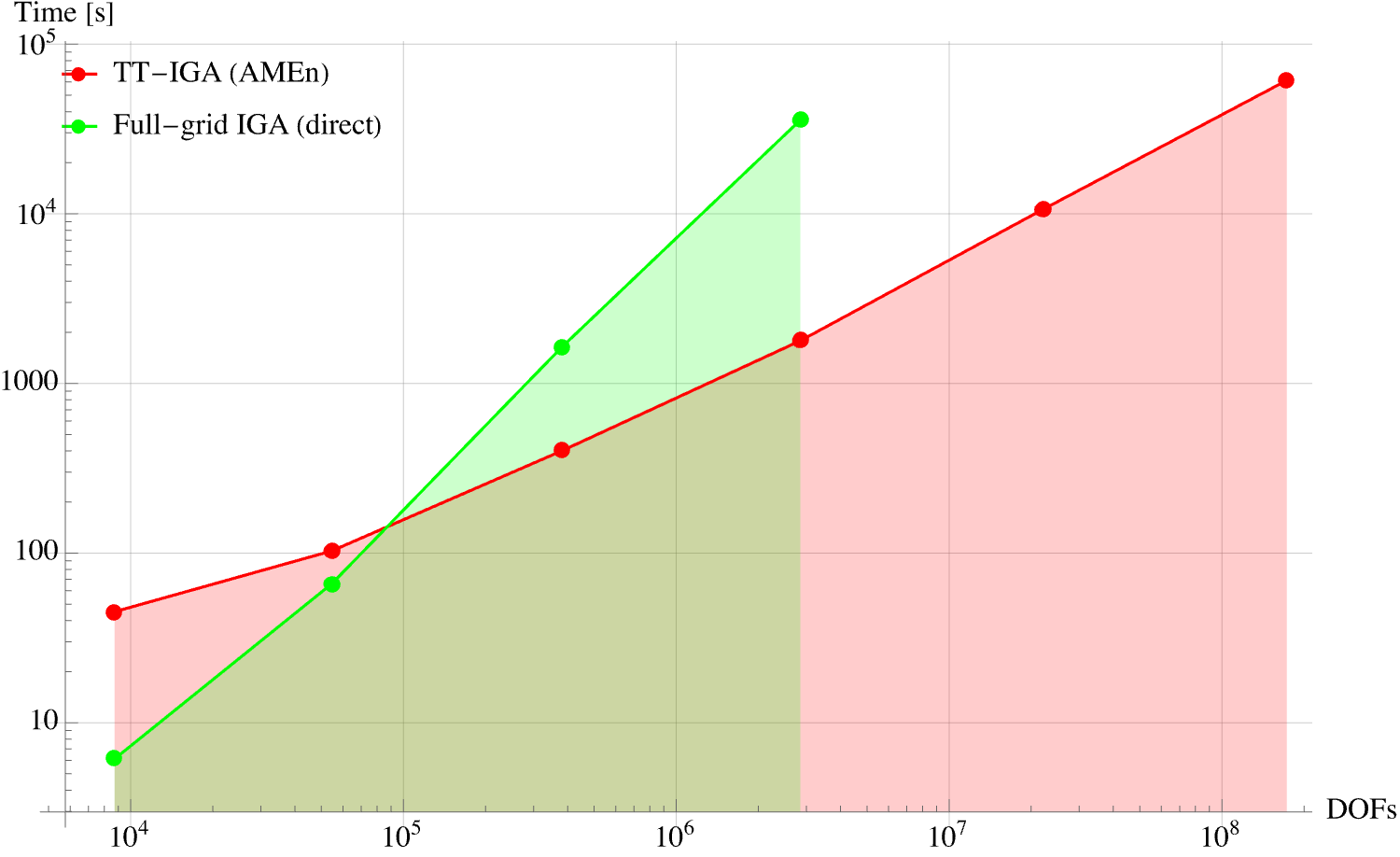}
\caption*{(d) Total computational time}
\endminipage
\caption{Compression ratios of the stiffness matrix $\boldsymbol{K}$, force vector $\boldsymbol{f}$, and displacement solution $\boldsymbol{u}$, together with TT-AMEn and full-grid computational times for the square-hole geometry.}
\label{Fig.square-hole-compression}
\end{center}
\end{figure}

\Figref{Fig.square-hole-compression} presents the compression ratios and solution times for this geometry. As in the open-hemisphere problem, the stiffness-operator compression increases rapidly with problem size, demonstrating that the geometry-dependent coupled elasticity operator remains strongly compressible in TT form. The displacement solution also becomes progressively more compressed under refinement, although its compression ratio is lower than that of the operator because the solution rank reflects the spatial complexity of the displacement field. The force vector remains inexpensive to represent because the prescribed traction is simple and is applied on a tensor-product boundary surface.

The solution-time comparison again reveals a crossover between the two approaches. The conventional sparse solver is competitive for the smallest discretizations, but its cost increases much faster as the system is refined. Once the problem reaches the intermediate-size regime, the TT-AMEn workflow becomes substantially faster and can continue to larger systems after the full-grid calculation is limited by the cost of sparse assembly and factorization. Together with the open-hemisphere results, this example shows that the favorable TT behavior is not restricted to the manufactured cylindrical test and persists for distinct curved three-dimensional geometries and loading conditions.
\section{Conclusion}
This work developed a direct tensor-train isogeometric formulation for coupled three-dimensional linear elasticity. Starting from the tensor-product spline representation of the geometry, the geometry-dependent elasticity coefficients are approximated by TT-cross interpolation and used to construct the coupled stiffness operator directly in TT format. The displacement components are incorporated through an additional TT mode, Dirichlet and Neumann boundary conditions are imposed without reverting to a full-grid representation, and the resulting elasticity system is solved in compressed form with the AMEn algorithm. The key innovation is the direct extension of fully assembled TT-IGA from a scalar PDE to a vector-valued three-dimensional elasticity system with coupled displacement components, while retaining the original spline geometry and spatial tensor-product structure. Tensor compression is therefore introduced at the level of the discretized elasticity formulation itself, rather than being applied only as a post-processing or reduced-order approximation to an already assembled full-grid system.

The numerical studies provide two main observations. First, the manufactured cylindrical example shows that the TT formulation reproduces the convergence behavior of the conventional full-grid IGA discretization. For quadratic basis functions, both approaches exhibit the same approximately third-order \(L_2\)-error convergence, demonstrating that the low-rank representation preserves the accuracy of the underlying discretization at the prescribed tolerances. Second, the storage advantage of the TT representation increases with problem size. Across the cylindrical, open-hemisphere, and square-hole geometries, the stiffness operator, force vector, and displacement field remain compressible while the corresponding full-grid storage grows rapidly. The solution-time comparisons show the complementary computational trend: TT-AMEn carries additional overhead for small systems, but becomes increasingly favorable under refinement as sparse full-grid assembly and direct solution become more expensive.

The large-scale examples highlight the practical consequence of this behavior. For the open hemispherical geometry, the TT formulation is carried to approximately \(1.6\times10^9\) degrees of freedom, well beyond the range for which the corresponding assembled full-grid calculations are practical in this study. Similar trends are observed for the square-hole geometry. These results indicate that, for tensor-product spline geometries that generate low-rank coefficient fields and solutions, the computational difficulty is governed less by the nominal full-grid dimension than by the TT ranks required to represent the coupled operator and displacement field. This is the central computational opportunity offered by the proposed formulation.

The method also has a clear scope. Its efficiency relies on tensor-product parameterizations and on sufficiently low TT ranks of the geometry-dependent coefficients, loading, boundary data, and solution. It is therefore not intended as a general replacement for arbitrary unstructured three-dimensional discretizations, and increasing geometric or material complexity may lead to rank growth that reduces the compression advantage. Future work should consequently focus on quantifying rank growth for increasingly nonseparable geometries and heterogeneous material fields, extending the formulation to multipatch and locally refined spline discretizations, and developing rank-adaptive strategies that balance approximation accuracy and computational cost. These extensions would broaden the class of structural problems for which a fully TT-based IGA formulation can provide an effective alternative to conventional full-grid simulation.
\appendix
\section{Computing storage}
\label{sec.Appen-Storage}
For all full-grid calculations reported in the numerical examples, stiffness storage is measured from the actual number of nonzero entries of the assembled sparse matrix. In the absence of trimming and \(G_0\) constraints, this quantity can be estimated analytically as
\begin{equation}
\operatorname{nnz}(\boldsymbol{K})
=
9\prod_{\alpha=1}^{3}
\left[n_\alpha(2p+1)-p(p+1)\right],
\end{equation}
where \(p\) is the spline degree and \(n_\alpha\) is the number of basis functions in direction \(\alpha\). The storage of a TT object is the total number of scalar entries in its cores, as defined in \Sref{sec.NotationTT}. For the stiffness operator, the compression ratio is therefore evaluated relative to the sparse full-grid storage,
\begin{equation}
\mathrm{CR}(\boldsymbol{K})
=
\frac{\operatorname{nnz}(\boldsymbol{K})}
{s_{\boldsymbol{K}^{\mathtt{TT}}}}.
\end{equation}
For a full-grid vector \(\boldsymbol{v}\) containing \(N\) scalar entries, such as the force or displacement vector, the corresponding compression ratio is
\begin{equation}
\mathrm{CR}(\boldsymbol{v})
=
\frac{N}{s_{\boldsymbol{v}^{\mathtt{TT}}}}.
\end{equation}
Thus, all compression ratios reported in the numerical examples compare the actual TT core storage with the storage of the corresponding conventional full-grid object.
\section{Manufactured solution for the 3D cylindrical ring}
\label{sec.Appen-Manufacture-Ring}
The strong form of the isotropic linear elasticity problem is
\begin{equation}
-\nabla\cdot \boldsymbol{\sigma}(\boldsymbol{u}) = \boldsymbol{f} \qquad \text{in } \Omega.
\end{equation}
Using the kinematic relation in \Eqref{eq.Kinematics} and the constitutive equation in \Eqref{eq.Constitutive}, with constant \(\lambda\) and \(\mu\), the strong form can be written as
\begin{equation}
\boldsymbol{f} = -\mu \Delta \boldsymbol{u} - (\lambda+\mu)\nabla(\nabla\cdot\boldsymbol{u}).    
\end{equation}
Consider a 3D cylindrical domain $\Omega = \left\{(x,y,z)\in\mathbb{R}^3: R_{in} \leq r=\sqrt{x^2+y^2} \leq R_{out},\quad 0 \leq z \leq t \right\}$. Let $R=x^2+y^2$ and choose the manufactured displacement as
\begin{equation}
\boldsymbol{u}_{\mathrm{ex}} = \begin{bmatrix}
u_x\\
u_y\\
u_z
\end{bmatrix}
= \begin{bmatrix}
\Phi\\
2\Phi\\
3\Phi
\end{bmatrix}
\end{equation}
where
\begin{equation}
\Phi(x,y,z) = P(x,y,z)S(x,y,z),
\end{equation}
with $P(x,y,z) = (R-r_i^2)(R-r_o^2)^2 z^2(t-z)^2$ and $S(x,y,z) = \sin(\alpha x)\cos(\beta y)\sin(\gamma z)$. The frequencies are chosen as $\alpha=\beta=\frac{\pi}{2r_o}$ and $\gamma=\frac{\pi}{t}$. On the boundaries $r = R_{in}$, $r = R_{out}$, $z = 0$, and $z = t$, we have
\begin{equation}
    \Phi=0 \qquad \nabla\Phi=\boldsymbol{0}
\end{equation}
which yields vanishing natural traction and zero displacement on the ring boundaries. Substituting the displacement field into the strong form gives the force vector as
\begin{equation}
\boldsymbol{f} = 
\begin{bmatrix}
&-\mu(Q_{xx}+Q_{yy}+Q_{zz})
-(\lambda+\mu)(Q_{xx}+2Q_{xy}+3Q_{xz})\\
&-2\mu(Q_{xx}+Q_{yy}+Q_{zz})
-(\lambda+\mu)(Q_{xy}+2Q_{yy}+3Q_{yz})\\
&-3\mu(Q_{xx}+Q_{yy}+Q_{zz})
-(\lambda+\mu)(Q_{xz}+2Q_{yz}+3Q_{zz}).
\end{bmatrix}
\end{equation}
where
\begin{equation}
A=R-r_i^2,\qquad
E=R-r_o^2,\qquad
B=E^2,\qquad
C=z^2(t-z)^2.
\end{equation}
and
\begin{equation}
\begin{split}
Q_{xx}
={}&CS\left[2E^2+16x^2E+A(4E+8x^2)-\alpha^2AE^2\right] +
2C\left[2xE^2+4xAE\right]
\alpha\cos(\alpha x)\cos(\beta y)\sin(\gamma z), \\
Q_{yy}
={}&
CS\left[
2E^2+16y^2E+A(4E+8y^2)-\beta^2AE^2
\right] - 2C\left[2yE^2+4yAE\right]
\beta\sin(\alpha x)\sin(\beta y)\sin(\gamma z), \\[0.5em]
Q_{zz}
={}&
AE^2\left[
C_{,zz}S
+2C_{,z}\gamma\sin(\alpha x)\cos(\beta y)\cos(\gamma z)
-\gamma^2CS
\right], \\[0.5em]
Q_{xy}
={}&
C\Big[
(16xyE+8Axy)S -(2xE^2+4xAE)\beta\sin(\alpha x)\sin(\beta y)\sin(\gamma z) \\
&+(2yE^2+4yAE)\alpha\cos(\alpha x)\cos(\beta y)\sin(\gamma z) -AE^2\alpha\beta\cos(\alpha x)\sin(\beta y)\sin(\gamma z)
\Big], \\[0.5em]
Q_{xz}
={}&
(2xE^2+4xAE)C_{,z}S +(2xE^2+4xAE)C\gamma\sin(\alpha x)\cos(\beta y)\cos(\gamma z) \\
&+AE^2C_{,z}\alpha\cos(\alpha x)\cos(\beta y)\sin(\gamma z)+AE^2C\alpha\gamma\cos(\alpha x)\cos(\beta y)\cos(\gamma z), \\[0.5em]
Q_{yz}
={}&
(2yE^2+4yAE)C_{,z}S +(2yE^2+4yAE)C\gamma\sin(\alpha x)\cos(\beta y)\cos(\gamma z) \\
&-AE^2C_{,z}\beta\sin(\alpha x)\sin(\beta y)\sin(\gamma z) -AE^2C\beta\gamma\sin(\alpha x)\sin(\beta y)\cos(\gamma z).
\end{split}
\end{equation}
The displacement error is measured using
\begin{equation}
\|\boldsymbol{u}_h-\boldsymbol{u}_{\mathrm{ex}}\|_{L^2(\Omega)}
=
\left[
\int_\Omega
\left(
(u_{h,x}-u_x)^2
+
(u_{h,y}-u_y)^2
+
(u_{h,z}-u_z)^2
\right)
\,d\Omega
\right]^{1/2}.
\end{equation}
The relative \(L^2\) error is
\begin{equation}
e_{L^2}^{\mathrm{rel}}
=
\frac{
\|\boldsymbol{u}_h-\boldsymbol{u}_{\mathrm{ex}}\|_{L^2(\Omega)}
}{
\|\boldsymbol{u}_{\mathrm{ex}}\|_{L^2(\Omega)}
}.
\end{equation}
\section{Application of boundary conditions}
\label{Sec.BCs}
\subsection{Imposition of Dirichlet boundary conditions}
Homogeneous Dirichlet conditions are imposed through a padding-mask construction that preserves the dimensions and tensor-product structure of the coupled TT system. The displacement components are stored together in one TT vector with ordering \(\{x,y,z,c\}\), where \(c=1,2,3\) denotes the displacement component. Three spatial TT masks are first constructed,
\[
\boldsymbol{S}_{x}^{s,\mathtt{TT}},\qquad
\boldsymbol{S}_{y}^{s,\mathtt{TT}},\qquad
\boldsymbol{S}_{z}^{s,\mathtt{TT}},
\]
with each mask obtained from the spatial TT identity by setting the prescribed boundary layers to zero in the corresponding row and column modes. Abstractly,
\[
\boldsymbol{S}_{a}^{s,\mathtt{TT}}
=
\mathcal{Z}\!\left(\boldsymbol{I}_{xyz}^{\mathtt{TT}};bc_a\right),
\qquad a\in\{x,y,z\},
\]
where \(\mathcal{Z}(\cdot;bc_a)\) denotes the operation of zeroing the layers associated with the boundary specification \(bc_a\).

The spatial masks are extended to the displacement-component mode using
\[
\boldsymbol{E}_{11}=
\begin{bmatrix}1&0&0\\0&0&0\\0&0&0\end{bmatrix},\qquad
\boldsymbol{E}_{22}=
\begin{bmatrix}0&0&0\\0&1&0\\0&0&0\end{bmatrix},\qquad
\boldsymbol{E}_{33}=
\begin{bmatrix}0&0&0\\0&0&0\\0&0&1\end{bmatrix}.
\]
The coupled padding mask is then
\begin{equation}
\boldsymbol{S}^{\mathtt{TT}}
=
\boldsymbol{S}_{x}^{s,\mathtt{TT}}\otimes\boldsymbol{E}_{11}
+
\boldsymbol{S}_{y}^{s,\mathtt{TT}}\otimes\boldsymbol{E}_{22}
+
\boldsymbol{S}_{z}^{s,\mathtt{TT}}\otimes\boldsymbol{E}_{33}.
\end{equation}
Using this mask, the constrained stiffness matrix and right-hand-side vector are constructed as
\begin{equation}
\boldsymbol{K}_{p}^{\mathtt{TT}}
=
\boldsymbol{S}^{\mathtt{TT}}
\boldsymbol{K}^{\mathtt{TT}}
\boldsymbol{S}^{\mathtt{TT}}
+
\alpha\left(\boldsymbol{I}^{\mathtt{TT}}-\boldsymbol{S}^{\mathtt{TT}}\right),
\end{equation}
and
\begin{equation}
\boldsymbol{f}_{p}^{\mathtt{TT}}
=
\boldsymbol{S}^{\mathtt{TT}}\boldsymbol{f}^{\mathtt{TT}}.
\end{equation}
Here, \(\boldsymbol{I}^{\mathtt{TT}}\) is the fourth-order TT identity matrix and \(\alpha\) is a prescribed diagonal value on constrained degrees of freedom; \(\alpha=10^{12}\) is used in the numerical examples. When a \(G_0\) continuity operation is present, the mask is constructed after that operation so that its mode sizes match those of the final TT operator. The constrained system solved by AMEn is therefore
\begin{equation}
\boldsymbol{K}_{p}^{\mathtt{TT}}
\boldsymbol{u}_{p}^{\mathtt{TT}}
=
\boldsymbol{f}_{p}^{\mathtt{TT}}.
\end{equation}

\subsection{TT construction of the Neumann force vector}
\label{App.TTNeumann}
In a three-dimensional problem, a Neumann boundary condition prescribed on a surface can be evaluated through a two-dimensional surface integral and then embedded into the full three-dimensional force vector. The same idea is used to construct the Neumann force vector in TT format. For a Neumann boundary condition prescribed on a surface \(\Gamma_N\), the corresponding force component is given by
\begin{equation}
f_i^N
=
\int_{\Gamma_N}
N_i\,\bar{t}\,d\Gamma,
\end{equation}
where \(N_i\) denotes the basis function and \(\bar{t}\) is the prescribed traction component. For the elasticity problem, the Neumann force has three components,
\begin{equation}
\boldsymbol{f}^{N,\mathtt{TT}}
=
\begin{bmatrix}
\boldsymbol{f}_{x}^{N,\mathtt{TT}}\\
\boldsymbol{f}_{y}^{N,\mathtt{TT}}\\
\boldsymbol{f}_{z}^{N,\mathtt{TT}}
\end{bmatrix}.
\end{equation}

For each Neumann boundary condition, we first select the corresponding two-dimensional boundary patch from the three-dimensional tensor-product grid. For example, if the boundary is the first layer in the \(\xi_3\) direction, the selected surface control points are
\begin{equation}
x_s(i,j)=x(i,j,1), \qquad
y_s(i,j)=y(i,j,1), \qquad
z_s(i,j)=z(i,j,1).
\end{equation}
On the selected surface, we define the weighted traction integrand
\begin{equation}
g(\xi_1,\xi_2)
=
\bar{t}
\bigl(
x(\xi_1,\xi_2),
y(\xi_1,\xi_2),
z(\xi_1,\xi_2)
\bigr)
J_s(\xi_1,\xi_2),
\end{equation}
where \(J_s\) is the surface Jacobian,
\begin{equation}
J_s
=
\left\|
\frac{\partial \boldsymbol{x}}{\partial \xi_1}
\times
\frac{\partial \boldsymbol{x}}{\partial \xi_2}
\right\|,
\end{equation}
and \(\partial \boldsymbol{x}/\partial \xi_1\) and \(\partial \boldsymbol{x}/\partial \xi_2\) are the covariant tangent vectors of the surface. The sampled function \(g\) is approximated in TT format as
\begin{equation}
g(\xi_1,\xi_2)
\approx
\sum_{\alpha=1}^{R}
G_g^{(1)}{}_{1,\xi_1,\alpha}
G_g^{(2)}{}_{\alpha,\xi_2,1}.
\end{equation}

For each TT rank term, we compute the one-dimensional force vectors
\begin{equation}
\boldsymbol{f}_{\alpha}^{(1)}
=
\int
\boldsymbol{N}^{(1)}(\xi_1)
G_g^{(1)}{}_{1,\xi_1,\alpha}
\,d\xi_1,
\end{equation}
and
\begin{equation}
\boldsymbol{f}_{\alpha}^{(2)}
=
\int
\boldsymbol{N}^{(2)}(\xi_2)
G_g^{(2)}{}_{\alpha,\xi_2,1}
\,d\xi_2.
\end{equation}
The surface force is then assembled as
\begin{equation}
\boldsymbol{f}^{s,\mathtt{TT}}
=
\sum_{\alpha=1}^{R}
\boldsymbol{f}_{\alpha}^{(1)}
\otimes
\boldsymbol{f}_{\alpha}^{(2)}.
\end{equation}
Equivalently, in TT notation,
\begin{equation}
\boldsymbol{f}^{s,\mathtt{TT}}(i,j)
=
\sum_{\alpha=1}^{R}
G_s^{(1)}{}_{1,i,\alpha}
G_s^{(2)}{}_{\alpha,j,1}.
\end{equation}

The scalar surface force is then assigned to the requested displacement component $(\boldsymbol{f}_{x}^{s,\mathtt{TT}}, \boldsymbol{f}_{y}^{s,\mathtt{TT}}, \boldsymbol{f}_{z}^{s,\mathtt{TT}})$. The surface force is finally embedded into the full three-dimensional tensor grid. For example, for the first boundary layer in the $\xi_3$ direction, the surface force has size \(n_1\times n_2\), and the corresponding volume force is defined by
\begin{equation}
\boldsymbol{f}^{v,\mathtt{TT}}(i,j,1)
=
\boldsymbol{f}^{s,\mathtt{TT}}(i,j),
\end{equation}
and
\begin{equation}
\boldsymbol{f}^{v,\mathtt{TT}}(i,j,k)
=
0,
\qquad
k\neq 1.
\end{equation}
The volume force is represented by
\begin{equation}
\boldsymbol{f}^{v,\mathtt{TT}}(i,j,k)
=
\sum_{\alpha=1}^{R}
G_v^{(1)}{}_{1,i,\alpha}
G_v^{(2)}{}_{\alpha,j,1}
G_v^{(3)}{}_{1,k,1}
\end{equation}
where
\begin{equation}
\begin{split}
&G_v^{(1)}{}_{1,i,\alpha} = G_s^{(1)}{}_{1,i,\alpha},\\
&G_v^{(2)}{}_{\alpha,j,1} = G_s^{(2)}{}_{\alpha,j,1},\\
&G_v^{(3)}{}_{1,k,1} = \delta_{k,1}.
\end{split}
\end{equation}
The same construction is used for the other directions and other boundary layers.
\section*{Acknowledgments}
The authors gratefully acknowledge the support of the Laboratory Directed Research and Development (LDRD) program of Los Alamos National Laboratory under project number 20260495ER. Los Alamos National Laboratory is operated by Triad National Security, LLC, for the National Nuclear Security Administration of U.S. Department of Energy (Contract No.\ 89233218CNA000001).

\section*{Data Availability}
No experimental datasets were used in this study. The numerical data supporting the findings of this work were generated from the simulations described in the manuscript and are available from the corresponding author upon reasonable request.

\bibliographystyle{IEEEtran}
\bibliography{References}
\end{document}